\documentclass{article}
\usepackage{amssymb,amsmath,
theorem,euscript}

\input{epsf}

\newcounter{sec}

\newcounter{punct}[sec]

\def\punct{\refstepcounter{punct}{\arabic{sec}.\arabic{punct}.  }}

\def\COUNTERS{\addtocounter{sec}{1}
              \setcounter{punct}{0}
          \setcounter{equation}{0}
          \setcounter{theorem}{0}
          \setcounter{problem}{0}
            \setcounter{Apunct}{0}
          }

\newtheorem{theorem}{Theorem}[sec]
\newtheorem{proposition}[theorem]{Proposition}

\newtheorem{lemma}[theorem]{Lemma}

\newtheorem{corollary}[theorem]{Corollary}
\newtheorem{observation}[theorem]{Observation}

\def\COUNTERS{\addtocounter{sec}{1}
              \setcounter{punct}{0}
          \setcounter{equation}{0}
          \setcounter{theorem}{0}
          }

\begin{document}

\def\SL{\mathrm {SL}}
\def\PSL{\mathrm {PSL}}
\def\SU{\mathrm {SU}}
\def\GL{\mathrm  {GL}}
\def\U{\mathrm  U}
\def\OO{\mathrm  O}
\def\Sp{\mathrm  {Sp}}
\def\SO{\mathrm  {SO}}
\def\SOS{\mathrm {SO}^*}

\def\PGL{\mathrm  {PGL}}
\def\PU{\mathrm {PU}}

\def\Gr{\mathrm{Gr}}

\def\Fl{\mathrm{Fl}}

\def\OSp{\mathrm {OSp}}

\def\Mat{\mathrm{Mat}}

\def\Pfaff{\mathrm {Pfaff}}

\def\Ind{\mathrm{Ind}}

\def\B{\mathbf B}

\def\phi{\varphi}
\def\epsilon{\varepsilon}
\def\kappa{\varkappa}

\def\le{\leqslant}
\def\ge{\geqslant}

\renewcommand{\Re}{\mathop{\rm Re}\nolimits}

\renewcommand{\Im}{\mathop{\rm Im}\nolimits}

\newcommand{\tr}{\mathop{\rm tr}\nolimits}

\newcommand{\im}{\mathop{\rm im}\nolimits}
\newcommand{\indef}{\mathop{\rm indef}\nolimits}
\newcommand{\dom}{\mathop{\rm dom}\nolimits}
\newcommand{\codim}{\mathop{\rm codim}\nolimits}

\def\cA{\mathcal A}
\def\cB{\mathcal B}
\def\cC{\mathcal C}
\def\cD{\mathcal D}
\def\cE{\mathcal E}
\def\cF{\mathcal F}
\def\cG{\mathcal G}
\def\cH{\mathcal H}
\def\cJ{\mathcal J}
\def\cI{\mathcal I}
\def\cK{\mathcal K}
\def\cL{\mathcal L}
\def\cM{\mathcal M}
\def\cN{\mathcal N}
\def\cO{\mathcal O}
\def\cP{\mathcal P}
\def\cQ{\mathcal Q}
\def\cR{\mathcal R}
\def\cS{\mathcal S}
\def\cT{\mathcal T}
\def\cU{\mathcal U}
\def\cV{\mathcal V}
\def\cW{\mathcal W}
\def\cX{\mathcal X}
\def\cY{\mathcal Y}
\def\cZ{\mathcal Z}

\def\frA{\mathfrak A}
\def\frB{\mathfrak B}
\def\frC{\mathfrak C}
\def\frD{\mathfrak D}
\def\frE{\mathfrak E}
\def\frF{\mathfrak F}
\def\frG{\mathfrak G}
\def\frH{\mathfrak H}
\def\frJ{\mathfrak J}
\def\frK{\mathfrak K}
\def\frL{\mathfrak L}
\def\frM{\mathfrak M}
\def\frN{\mathfrak N}
\def\frO{\mathfrak O}
\def\frP{\mathfrak P}
\def\frQ{\mathfrak Q}
\def\frR{\mathfrak R}
\def\frS{\mathfrak S}
\def\frT{\mathfrak T}
\def\frU{\mathfrak U}
\def\frV{\mathfrak V}
\def\frW{\mathfrak W}
\def\frX{\mathfrak X}
\def\frY{\mathfrak Y}
\def\frZ{\mathfrak Z}

\def\fra{\mathfrak a}
\def\frb{\mathfrak b}
\def\frc{\mathfrak c}
\def\frd{\mathfrak d}
\def\fre{\mathfrak e}
\def\frf{\mathfrak f}
\def\frg{\mathfrak g}
\def\frh{\mathfrak h}
\def\fri{\mathfrak i}
\def\frj{\mathfrak j}
\def\frk{\mathfrak k}
\def\frl{\mathfrak l}
\def\frm{\mathfrak m}
\def\frn{\mathfrak n}
\def\fro{\mathfrak o}
\def\frp{\mathfrak p}
\def\frq{\mathfrak q}
\def\frr{\mathfrak r}
\def\frs{\mathfrak s}
\def\frt{\mathfrak t}
\def\fru{\mathfrak u}
\def\frv{\mathfrak v}
\def\frw{\mathfrak w}
\def\frx{\mathfrak x}
\def\fry{\mathfrak y}
\def\frz{\mathfrak z}

\def\bfa{\mathbf a}
\def\bfb{\mathbf b}
\def\bfc{\mathbf c}
\def\bfd{\mathbf d}
\def\bfe{\mathbf e}
\def\bff{\mathbf f}
\def\bfg{\mathbf g}
\def\bfh{\mathbf h}
\def\bfi{\mathbf i}
\def\bfj{\mathbf j}
\def\bfk{\mathbf k}
\def\bfl{\mathbf l}
\def\bfm{\mathbf m}
\def\bfn{\mathbf n}
\def\bfo{\mathbf o}
\def\bfp{\mathbf q}
\def\bfr{\mathbf r}
\def\bfs{\mathbf s}
\def\bft{\mathbf t}
\def\bfu{\mathbf u}
\def\bfv{\mathbf v}
\def\bfw{\mathbf w}
\def\bfx{\mathbf x}
\def\bfy{\mathbf y}
\def\bfz{\mathbf z}

\def\bfA{\mathbf A}
\def\bfB{\mathbf B}
\def\bfC{\mathbf C}
\def\bfD{\mathbf D}
\def\bfE{\mathbf E}
\def\bfF{\mathbf F}
\def\bfG{\mathbf G}
\def\bfH{\mathbf H}
\def\bfI{\mathbf I}
\def\bfJ{\mathbf J}
\def\bfK{\mathbf K}
\def\bfL{\mathbf L}
\def\bfM{\mathbf M}
\def\bfN{\mathbf N}
\def\bfO{\mathbf O}
\def\bfP{\mathbf P}
\def\bfQ{\mathbf Q}
\def\bfR{\mathbf R}
\def\bfS{\mathbf S}
\def\bfT{\mathbf T}
\def\bfU{\mathbf U}
\def\bfV{\mathbf V}
\def\bfW{\mathbf W}
\def\bfX{\mathbf X}
\def\bfY{\mathbf Y}
\def\bfZ{\mathbf Z}

\def\R {{\mathbb R }}
 \def\C {{\mathbb C }}
  \def\Z{{\mathbb Z}}
  \def\H{{\mathbb H}}
\def\K{{\mathbb K}}
\def\N{{\mathbb N}}
\def\Q{{\mathbb Q}}
\def\A{{\mathbb A}}

\def\T{\mathbb T}

\def\bbA{\mathbb A}
\def\bbB{\mathbb B}
\def\bbD{\mathbb D}
\def\bbE{\mathbb E}
\def\bbF{\mathbb F}
\def\bbG{\mathbb G}
\def\bbI{\mathbb I}
\def\bbJ{\mathbb J}
\def\bbL{\mathbb L}
\def\bbM{\mathbb M}
\def\bbN{\mathbb N}
\def\bbO{\mathbb O}
\def\bbP{\mathbb P}
\def\bbQ{\mathbb Q}
\def\bbS{\mathbb S}
\def\bbT{\mathbb T}
\def\bbU{\mathbb U}
\def\bbV{\mathbb V}
\def\bbW{\mathbb W}
\def\bbX{\mathbb X}
\def\bbY{\mathbb Y}

\def\la{\langle}
\def\ra{\rangle}

 \def\ov{\overline}
\def\wt{\widetilde}
\def\wh{\widehat}

\def\P{\mathbb P}

\def\bO{\bf O}

\def\arr{\rightrightarrows}

\def\SS{\smallskip}
\def\MS{\medskip}

\def\ev{{\mathrm{even}}}
\def\od{{\mathrm{odd}}}

\def\q{\quad}

\def\F{\mathbf F}

\def\b{\mathbf b}

\def\RA{\Longrightarrow}

\def\BS{\bigskip}

\def\bm{{\mathbf m}}
\def\bl{{\mathbf l}}

\def\fD{{\frak D}}
\def\frk{\mathfrak k}

\def\B{{\rm B}}

\def\b{C}
\def\T{T}

\def\phi{\varphi}
\def\epsilon{\varepsilon}
\def\kappa{\varkappa}
\def\le{\leqslant}
\def\ge{\geqslant}

\def\konets{\hfill$\boxtimes$}

\newcommand{\Ker}{\mathop{\rm Ker}\nolimits}
  \renewcommand{\Im}{\mathop{\rm Im}\nolimits}
    \newcommand{\Dom}{\mathop{\rm Dom}\nolimits}
     \newcommand{\Indef}{\mathop{\rm Indef}\nolimits}
\newcommand{\graph}{\mathop{\rm graph}\nolimits}

\def\const{{\rm const}}

\def\tto{\rightrightarrows}

\renewcommand{\Re}{\mathop{\rm Re}\nolimits}
\renewcommand{\Im}{\mathop{\rm Im}\nolimits}

\def\palka{\vphantom{1^G}}



\def\kartinka{

\pagebreak

\unitlength=1mm

\begin{picture}(100,100)
\put(10,50){\line(0,-1){43}}
\put(20,40){\line(0,-1){33}}
\put(30,30){\line(0,-1){23}}
\put(40,20){\line(0,-1){13}}
\put(50,10){\line(0,-1){3}}

\put(10,50){\line(-1,0){3}}
\put(20,40){\line(-1,0){13}}
\put(30,30){\line(-1,0){23}}
\put(40,20){\line(-1,0){33}}
\put(50,10){\line(-1,0){43}}

\multiput(10,80)(10,-10){8}{\circle*{1}}
\multiput(20,80)(10,-10){8}{\circle*{1}}
\multiput(30,80)(10,-10){8}{\circle*{1}}
\multiput(40,80)(10,-10){8}{\circle*{1}}

\put(60,60){\vector(0,1){30}}
\put(60,60){\vector(1,0){30}}
\put(90,65){$\sigma$}
\put(65,90){$\tau$}

{\linethickness{1mm}
\multiput(10,50)(10,-10){5}
{\line(0,1){10}}
\multiput(20,50)(10,-10){5}%
{\line(0,1){10}}
\multiput(10,50)(10,-10) {5}%
{\line(1,0){10}}
\multiput(10,60)(10,-10) {5}%
{\line(1,0){10}}
\put(10,60){\line(0,1){10}}
\put(60,10){\line(1,0){10}}
\put(60,10){\line(0,-1){3}}
\put(70,10){\line(0,-1){3}}
\put(10,60){\line(-1,0){3}}
\put(10,70){\line(-1,0){3}}
}

\multiput(11,51)(10,-10){5}
{%
\multiput(0,0)(1,0){9}%
{\multiput(0,0)(0,1){9}%
{\circle*{0.2}}%
}%
}

\put(7,63){\line(1,-1){56}}

\put(70,60){\line(0,1){2}}
\put(70,63){1}
\put(60,70){\line(1,0){2}}
\put(63,70){1}

\put(35,35){\circle*{2}}

\put(100,30){\vector(1,-1){10}}
\put(105,30){Shift}

\multiput(0,10)(0,10){8}\dots
\multiput(10,0)(10,0){11}\vdots

\end{picture}

{\sf
 1. The dotted squares correspond to unitary representations
$\rho_{\sigma|\tau}$.

2. Vertical and horizontal rays in the south-west of Figure
correspond to nondegenerate
 highest weight and lowest weight representations.
Fat points  correspond
to degenerated highest and lowest weight
  representations, and also to the unipotent representations.
  The point
$(\sigma,\tau)=(0,0)$ corresponds
to the trivial one-dimensional representation.

3. In points of the thick segments,
we have some exotic unitary sub-quotients.

4. The shift
$(\sigma,\tau)\mapsto(\sigma+1,\tau-1)$
 send a  representation
$\rho_{\sigma|\tau}$ of $\SU(n,n)^\sim$
to an equivalent representation.

5. The permutation of the axes
$(\tau,\sigma)\mapsto(\sigma,\tau)$
 gives a complex conjugate representation.

6. The symmetry with respect to the point $(-n/2,-n/2)$
(black circle) gives a dual representation
(for odd $n$ this point is a center of a dotted square;
for even $n$ this point is  a common vertex of two
dotted squares).

7. For $\sigma+\tau=n$ (the  diagonal line)
our Hermitian
form is the standard $L^2$-product.

8. Linear (non-projective) representations
of $\U(n,n)$ correspond to  the
family of parallel lines
$\sigma-\tau\in\Z$.}
\pagebreak
}

\begin{center}

{\Large \bf

Stein--Sahi complementary series\\
 and their degenerations}

\vspace{22pt}

{\large\sc

Yuri A. Neretin%
\footnote{Supported by the grant FWF, project P19064,
 Russian Federal Agency for Nuclear Energy,
Dutch  grant NWO.047.017.015, and grant JSPS-RFBR-07.01.91209}
}

\end{center}

{\small The paper is an introduction
to the Stein--Sahi complementary series,
 and the unipotent representations.
We also discuss some open problems related
to these objects.
For the sake of simplicity,
we consider only the groups $\U(n,n)$.}

\section{Introduction}

\COUNTERS

 This paper%
\footnote{It is a strongly revised version
of two sections of my preprint \cite{Ner-preprint}.}
 is an attempt to
present  an introduction to the Stein-Sahi complementary
series  available for non-experts and beginners.


\SS

{\bf\punct History of the subject.}
 Theory of infinite dimensional
representations of semi-simple groups
 was initiated in pioneer works
of I.~M.~Gelfand and M.~A.~Naimark (1946--1950),
 V.~Bargmann \cite{Bar} (1947), and
 K.~O.~Friedrichs \cite{Fri} (1951--1953).
The book \cite{GN2} by I.~M.~Gelfand and M.~A.~Naimark
(1950) contains a  well-developed theory
for  complex classical
 groups $\GL(n,\C)$, $\SO(n,\C)$, $\Sp(2n,\C)$
(the parabolic induction, complementary series,
 spherical functions, characters, Plancherel theorems).
However, this classical book%
\footnote{Unfortunately the book  exists only in Russian
and German.} contained various statements
and asseverations that were not actually proved.
In the modern terminology, some of chapters were
'mathematical physics'.
The most of these  statements were really proved by 1958--1962
in works of different authors
(Harish-Chandra, F.~A.~Berezin, etc.).

In particular, I.~M.~Gelfand and M.~A.~Naimark (1950)
 claimed that they classified all unitary representations
of $\GL(n,\C)$, $\SO(n,\C)$, $\Sp(2n,\C)$.
 E.~Stein \cite{Ste} compared Gelfand--Naimark
constructions  for groups $\SL(4,\C)\simeq\SO(6,\C)$
and observed that they are not equivalent.
In 1967 E.~Stein constructed 'new'
unitary representations of $\SL(2n,\C)$.

D.~Vogan \cite{Vog} in 1986
obtained the classification
of unitary representations
 of groups $\GL(2n)$ over real numbers $\R$
and quaternions $\H$. In particular, this work contains
extension of  Stein's construction to
these groups.
In 1990s, the Stein-type representations were a
 topic of interest  of S.~Sahi
 see \cite{Sah1},
 \cite{Sah2},  \cite{Sah3}, S.~Sahi--E.~Stein
 \cite{SS}, A.~Dvorsky--S.~Sahi \cite{DS1}--\cite{DS2}.
 In particular, Sahi
extended the construction to other series
of classical groups, precisely to
the groups
$\SO(2n,2n)$, $\U(n,n)$, $\Sp(n,n)$,
$\Sp(2n,\R)$, $\SOS(4n)$, $\Sp(4n,\C)$, and
$\SO(2n,\C)$.


\SS

{\bf\punct Stein--Sahi representations for $\U(n,n)$.}
Denote by $\U(n)$ the group of unitary $n\times n$-matrices.
Consider the {\it pseudo-unitary} group $\U(n,n)$.
 We realize it as the group
of $(n+n)\times (n+n)$-matrices
$g=\begin{pmatrix}a&b\\c&d\end{pmatrix}$
satisfying the condition
$$
\begin{pmatrix}a&b\\c&d\end{pmatrix}
\begin{pmatrix}1&0\\0&-1\end{pmatrix}
\begin{pmatrix}a&b\\c&d\end{pmatrix}^*
=
\begin{pmatrix}1&0\\0&-1\end{pmatrix}
.
$$

\begin{lemma} The  formula
\begin{equation}
z\mapsto z^{[g]}:=(a+zc)^{-1} (b+zd)
\label{eq:action-un}
\end{equation}
determines an action of the group
$\U(n,n)$ on the space $\U(n)$.
\end{lemma}

The unitary group is equipped by the Haar measure
$d\mu(z)$, hence we can determine the Jacobian of a transformation
(\ref{eq:action-un}) by
$$
J(g,z)=\frac{d\mu(z^{[g]})}{d\mu(z)}
.
$$

\begin{lemma}
 The Jacobian of the transformation
$z\mapsto z^{[g]}$ on $\U(n)$ is given by
$$
J(g,z)=
|\det(a+zc)|^{-2n}
.$$
\end{lemma}

Fix $\sigma$, $\tau\in\C$.
For $g\in\U(n,n)$ we define the following linear operator
in the space $C^\infty(\U(n))$:
\begin{equation}
\rho_{\sigma|\tau}(g) f(z)=
f(z^{[g]})
\det(a+zc)^{-n-\tau}
\det\ov{(a+zc)}^{-n-\sigma}
.\end{equation}

The formula includes powers of complex numbers,
precise definition is given below.
In fact, $g\mapsto\rho_{\sigma|\tau}(g)$
 is a well-defined operator-valued function
on the universal covering
group $\U(n,n)^\sim$ of $\U(n,n)$.

The chain rule for Jacobians,
\begin{equation}
J(g_1g_2,z)=J(g_1,z)J(g_2,z^{[g_1]})
\label{eq:chain-rule},
\end{equation}
implies
$$
\rho_{\sigma|\tau}(g_1)\rho_{\sigma|\tau}(g_2)
=\rho_{\sigma|\tau}(g_1g_2)
.
$$
In other words, $\rho_{\sigma|\tau}$ is
 a linear representation of the
group $\U(n,n)^\sim$.


\begin{observation}
\label{obs:1.2}
If $\Re \sigma+\Re \tau=-n$, $\Im\sigma=\Im\tau$
then a representation $\rho_{\sigma|\tau}$
is unitary in $L^2(\U(n))$.
\end{observation}

This easily follows from the formula for the Jacobian.

\SS

Next, let $\sigma$, $\tau$ be real.
We define the Hermitian form on $C^\infty(\U(n))$
by the formula
\begin{equation}
\la f_1,f_2\ra_{\sigma|\tau}:
=
\int_{\U(n)}\int_{\U(n)}
\det(1-zu^*)^\sigma (1-z^*u)^\tau
f_1(z)\,\ov{f_2(u)}\, d\mu(z)\,d\mu(u)
.\end{equation}

\begin{proposition}
The operators $\rho_{\sigma|\tau}(g)$ preserve
the Hermitian form $\la \cdot,\cdot\ra_{\sigma|\tau}$.
\end{proposition}

\begin{theorem}
For $\sigma$, $\tau\not\in \Z$,
the Hermitian form $\la \cdot,\cdot\ra_{\sigma|\tau}$
is positive iff integer parts
of numbers $-\sigma-n$ and $\tau$ are equal.
\end{theorem}

In fact, the domain of positivity is
 the square $-1<\tau<0$, $-n<\sigma<-n+1$
and its shifts by vectors $(-j,j)$, $j\in \Z$
, see Figure \ref{fig:big}.

In particular, under this condition,
a representation $\rho_{\sigma|\tau}$ is unitary.

\SS

For some values of $(\sigma,\tau)$ the form
$\la \cdot,\cdot\ra_{\sigma|\tau}$
is positive semi-definite. There are two the most important
cases.

\SS

1. For $\tau=0$, we get
 highest weight representations
(or holomorphic representations). Thus, the Stein--Sahi
representations are nearest relatives of
holomorphic representations.

\SS

2. For $\tau=0$, $\sigma=0$, $-1$, $-2$, \dots, $-n$
we obtain some exotic 'small' representations of $\U(n,n)$.

\SS


{\bf\punct The structure of the paper.}
{\it We discuss only groups%
\footnote{A comment for experts.
Stein--Sahi representations of a semisimple
Lie group $G$ are complementary
series induced
from a maximal parabolic subgroup with Abelian nilpotent radical.
\newline
$\vphantom{bb}$ \qquad
 The cases $G=\U(n,n)$,
 $\Sp(2n,\R)$,
$G=\SOS(4n)$ (related to tube type Hermitian symmetric spaces)
 are parallel. The  only difficulty is Theorem \ref{th:main}
(the expansion of the integral kernel in characters, we choose
$G=\U(n,n)$, because this can be done by  elementary tools).
In the general Hermitian case, one can refer to the version
of the Kadell integral \cite{Kad} from \cite{Ner-stein}
(the integrand is
a product of a Jack polynomial and a Selberg-type factor.
\newline
$\vphantom{bb}$ \qquad
 For other series of groups,
 Stein-Sahi representations
depend on one parameter, and picture is more pure
(in particular, inner products for degenerate ('unipotent')
 representations can be written immediately).
A $BC$-analog of Kadell integral is unknown
(certainly, it must exist, and some special cases
were evaluated in the literature, see e.g.,\cite{Ner-preprint}).
On the other hand, Stein-Sahi
representations have multiplicity free $K$-spectra. In such
situation, there is lot of ways for examination of positivity
of inner products, see e.g. \cite{Sah2}, \cite{Sah3},
\cite{BOO}.
\newline
$\vphantom{bb}$ \qquad
New elements of this paper are 'blow-up construction' for unipotent
representations and (apparently) tame models for representations
of universal coverings. The representations themselves
were constructed in works of  Sahi.}
 $\U(n,n)$.}

In Section 2 we consider the case $n=1$
and present the
Pukanszky classification \cite{Puk-classic}
 of unitary representations
of the universal covering group of
$\SL(2,\R)\simeq\SU(1,1)$.

In Section 3 we discuss Stein-Sahi representations
of arbitrary $\U(n,n)$.
In Section 4 we explain relations of Stein--Sahi representations
and holomorphic representations.
In Section 5 we give explicit constructions of
the Sahi 'unipotent' representations.

In Section 6 we discuss some open
 problems of harmonic analysis.

\SS


{\bf\punct Notation.} Let $a$, $u$, $v\in\C$.
Denote
\begin{equation}
a^{\{u|v\}}:=a^u\ov a^v
\label{eq:double-power}
.
\end{equation}
If $u-v\in\Z$, then this expression
 is well defined for all $a\ne 0$.
However, the expression is well defined in many other
situations, for instance
if $|1-a|<1$ and  $u$, $v$ are arbitrary
(and even for $|1-a|= 1$, $a\ne 1$)

\SS

The {\it norm} $\|z\|$ of an $n\times n$-matrix $z$
is the usual norm
of a linear operator in the standard Euclidean space $\C^n$.

\SS

We denote the Haar measure on the unitary group $\U(n)$
by $\mu$; assume that the complete measure of the group is 1.

\SS

The {\it Pochhammer symbol}  is given by
\begin{equation}
(a)_n:=\frac{\Gamma(a+n)}{\Gamma(a)}
=\begin{cases}
a(a+1)\dots(a+n-1)\qquad&\text{if $n\ge 0$}
\\
\frac{1}{(a-1)\dots(a-n)}\qquad&\text{if $n< 0$}
.
\end{cases}
\label{eq:pochhammer}
\end{equation}

\section{Unitary representations
 of $\SU(1,1)$}

\COUNTERS

Denote by
 $\SU(1,1)^\sim$
the universal covering  group of $\SU(1,1)$.

In this section, we present constructions of
all irreducible unitary representations
of  $\SU(1,1)^\sim$.
According the Bargmann--Pukanszky theorem there are
4 types of such representations:

\SS

a) unitary principal series;

\SS

b) complementary series;

\SS

c) highest weight and lowest weight
representations;

\SS

d) The one-dimensional representation.

\SS

Models of these representations are given  below.

\SS

The general Stein--Sahi representations are a strange
'higher copy' of the $\SU(1,1)$-picture.

\SS

{\bf References.} The classification of unitary representations
of $\SL(2,\R)\simeq\SU(1,1)$
 was obtained by V.~Bargmann \cite{Bar};
it was extended to the
 $\SU(1,1)^\sim$ by L.~Pukanszky \cite{Puk-classic},
 see also P.~Sally \cite{Sal}. \hfill $\square$


\BS

{\bf\large A. Preliminaries}

\BS


{\bf \punct Fourier series and distributions.}
By $S^1$ we denote  the unit circle
$|z|=1$ in the complex plane $\C$.
We parameterize  $S^1$ by   $z=e^{i\phi}$.

By $C^\infty(S^1)$ we denote  the space
 of smooth functions on $S^1$. Recall, that
 $$
f(\phi)=\sum_{n=-\infty}^\infty
a_n e^{in\phi}\in C^\infty(S^1)\qquad \text{iff
$|a_n|=o(|n|^{-L})$ for all $L$.}
 $$

Recall that a distribution $h(\phi)$ on the circle
admits an expansion into a Fourier series,
$$
h(\phi)=\sum_{n=0}^\infty b_n e^{in\phi}
,
\quad
\text{where $|b_n|=O(|n|^L)$ for some $L$.}
$$

For $s\in\R$ we define the
 {\it Sobolev space} $W^s(S^1)$ as
the space of distributions
$$
h(\phi)=\sum_{n=0}^\infty b_n e^{in\phi}
\quad
\text{such that $\sum |b_n|^2 (1+|n|)^{2s}<\infty$.}
$$
By definition, $W^0(S^1)=L^2(S^1)$.
For positive integer $s=k$
this condition is equivalent
$\frac{\partial^k}{\partial \phi^k} h\in L^2(S^1)$.
Evidently, $s<s'$ implies $W^s\supset W^{s'}$.

\SS


{\bf \punct The group  $\SU(1,1)$.}
The group $\SU(1,1)\simeq\SL(2,\R)$
consists of all complex $2\times 2$-matrices
having the form
$$
g=\begin{pmatrix}a&b\\ \ov b&\ov a\end{pmatrix}
,\qquad \text{where $|a|^2-|b|^2=1$.}
$$
This group acts on the disc $|z|<1$ and on the circle
$|z|=1$ by the {\it M\"obius transformations}
$$
z\mapsto (a+\ov b z)^{-1}(b+\ov a z)
.
$$


{\bf \punct A model of the universal covering group $\SU(1,1)^\sim$.}
Recall that the fundamental group of
$\SU(1,1)$ is $\Z$.
A loop generating the fundamental group
is
\begin{equation}
\frR(\phi)=
\begin{pmatrix}e^{i\phi}&0\\0&e^{-i\phi}\end{pmatrix}
,
\qquad \frR(2\pi)=\frR(0)=1.
\label{eq:R-phi}
\end{equation}
Some example of  multi-valued continuous function
on $\SU(1,1)$ are
$$
\begin{pmatrix}a&b\\ \ov b&\ov a\end{pmatrix}
\mapsto \ln a
,\qquad
\begin{pmatrix}a&b\\ \ov b&\ov a\end{pmatrix}
\mapsto
a^\lambda:=a^{\lambda\ln a}
.
$$


We can realize $\SU(1,1)^\sim$ as a subset
in $\SU(1,1)\times \C$ consisting of pairs
$$
\left(\,\begin{pmatrix}a&b\\\ov b&\ov a\end{pmatrix},
\, \sigma\right),
\quad\text{where $e^\sigma=a$.}
$$
Thus, for a given matrix
$\begin{pmatrix}a&b\\\ov b&\ov a\end{pmatrix}$
the parameter $\sigma$ ranges if the countable set
$\sigma=\ln a+ 2\pi ki$.

Define a multiplication in $\SU(1,1)\times \C$
by
$$
(g_1,\sigma_1)\circ(g_2,\sigma_2)=
(g_1g_2,\sigma_1+\sigma_2+c(g_1,g_2))
,$$
where $c(g_1,g_2)$ is the {\it Berezin--Guichardet} cocycle,
$$
c(g_1, g_2)=\ln\frac{a_3}{a_1 a_2}
,
$$
Here $a_3$ is the matrix element of $g_3=g_1g_2$.

\begin{theorem}
a) $\left|\frac{a_3}{a_1a_2}-1\right|<1$,
 and therefore the logarithm is well defined.

\SS

b) The operation $\circ$ determines the structure
of a group on $\SU(1,1)\times\C$.

\SS

c) $\SU(1,1)^\sim$ is a subgroup  in the latter group.
\end{theorem}

The proof is a simple and nice exercise.

\SS

Now we can define the single-valued function $\ln a$
on $\SU(1,1)$ by setting
$\ln a:=\sigma$.


\BS

{\bf\large B. Non-unitary and unitary principal series}

\BS


{\bf\punct Principal series
of representations of $\SU(1,1)$.}
Fix $p$, $q\in\C$. For $g\in\SU(1,1)$ define the operator
$T_{p|q}(g)$
in the space $C^\infty(S^1)$ by the formula
\begin{equation}
T_{p|q}
\begin{pmatrix}a&b\\ \ov b&\ov a\end{pmatrix}
f(z)
=f\Bigl(\frac{b+\ov a z}{a+\ov b z}\Bigr)
(a+\ov b z)^{\{-p|-q\}}
\label{eq:tpq}
,\end{equation}
here we use the notation (\ref{eq:double-power})
for complex powers.

\begin{observation}
a) $T_{p|q}$ is a well-defined operator-valued
function on $\SU(1,1)^\sim$.

\SS

b) It satisfies
$$
T_{p|q}(g_1)T_{p|q}(g_2)=T_{p|q}(g_1g_2)
.$$
\end{observation}

{\sc Proof.} a)  First,
$$
(a+\ov b z)^{-p} \ov{(a+\ov b z)}^{\,\,-q}=
a^{-p}\cdot \ov a^{\,\,-q}(1+a^{-1}\ov b z)^{-p}
\ov{(1+a^{-1}\ov b z)}^{\,\,-q}
.$$
Since $|z|=1$ and $|a|>|b|$, the last two factors
are well defined. Next,
$$
a^{-p}\, \ov a^{\,\,-q}:=\exp\Bigl\{- p\ln a +q\,\ov{\ln a}\Bigr\}
$$
and $\ln a$ is a well-defined function on
$\SU(1,1)^\sim$.

Proof of b).
One can verify this identity for $g_1$, $g_2$
near the unit and refer to the analytic continuation.
\hfill $\square$

\SS

The representations $T_{p|q}(g)$ are called
{\it  representations of the principal (non-unitary)
series}.

\SS

{\sc Remark.} a) A representation $T_{p|q}$ is a
 single-valued representation
of $\SU(1,1)$ iff $p-q$ is integer.


\SS

{\bf\punct  The action of the Lie algebra.}
The Lie algebra $\frs\fru(1,1)$
of $\SU(1,1)$ consists of matrices
$$
\begin{pmatrix}
i\alpha&\beta\\\ov\beta&-i\alpha
\end{pmatrix},
\qquad
\text{where $\alpha\in\R$, $\beta\in\C$.}
$$
It is convenient to take the following basis
 in the complexification $\frs\fru(1,1)_\C=\frs\frl(2,\C)$:
\begin{equation}
L_0:=
\frac 12\begin{pmatrix}-1&0\\0&1\end{pmatrix},
\quad
L_-:=\begin{pmatrix} 0&1\\0&0\end{pmatrix},
\quad
L_+:=\begin{pmatrix} 0&0\\-1&0\end{pmatrix}
\label{eq:LLL1}
\end{equation}

These generators
act in $C^\infty(S^1)$ by the following operators
\begin{equation}
L_0=z\frac d{dz}+\frac 12(p-q),
\qquad
L_-=\frac d{dz}-qz^{-1},
\qquad
L_+=z^2\frac d{dz}+pz.
\label{eq:LLL2}
\end{equation}
Equivalently,
\begin{equation}
L_0 z^n=\bigl(n+\frac12(p-q)\bigr)z^n,
\qquad
L_- z^n=(n-q)z^{n-1},
\qquad
L_+ z^n=(n+p) z^{n+1}
.
\end{equation}


{\bf\punct Subrepresentations.}

\begin{proposition}
A representation $T_{p|q}$ is irreducible iff
$p$, $q\notin \Z$.
\end{proposition}

{\sc Proof.} Let $p$, $q\notin\Z$.  Consider an
$L_0$-eigenvector $z^n$. Then
all vectors $(L_+)^k z^n$, $(L_-)^l z^n$
are nonzero. They span the whole space $C^\infty(S^1)$.
\hfill $\square$

\begin{observation}
a) If $q\in \Z$, then $z^q$, $z^{q+1}$, \dots span
a subrepresentation in $T_{p|q}$.

\SS

b) If $p\in \Z$, then $z^{-p}$, $z^{-p-1}$, $z^{-p-2}$,
\dots
span a subrepresentation in $T_{p|q}$.
\end{observation}

{\sc Proof of a).} Clearly, our subspace is
$L^0$-invariant and $L^+$-invariant. On the other hand,
$L^-z^q=0$, and we can not leave our subspace.%
\hfill $\square$

\SS

\begin{figure}
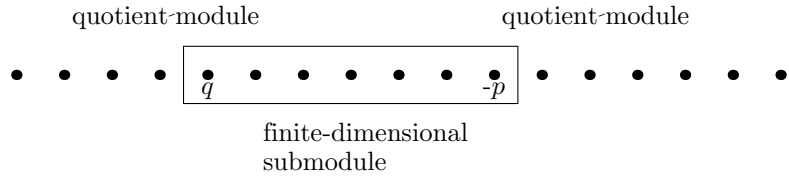


\epsfbox{sobolev.1}

a) {\it $q$ is integer};

\bigskip

\epsfbox{sobolev.2}

b) {\it $p$ is integer};

\bigskip

\epsfbox{sobolev.3}

c) {\it $p$, $q$ are integer,  $q+p\ge 1$};

\bigskip

\epsfbox{sobolev.4}

d) {\it $p$, $q$ are integer,  $q+p\le 1$}.

\caption{Subrepresentations of the principal series.
Black circles enumerate vectors $z^n$. A representation
$T_{p|q}$
is reducible iff $p\in\Z$ or $q\in\Z$.
\label{fig:su11-subr}}
\end{figure}

All  possible positions of subrepresentations
of $T_{p|q}$ are listed on Figure \ref{fig:su11-subr}.


\SS

{\bf\punct Shifts of parameters.}

\begin{observation}
\label{obs:u11-shift}
If $k$ is integer, then
$T_{p+k|q-k}\simeq T_{p|q}$.
The intertwining operator is
$$
Af(z)=z^k f(z)
.
$$
\end{observation}

A verification is trivial.
\hfill $\square$


\SS

{\bf \punct Duality.}
Consider the bilinear  map
$$
\Pi:C^\infty(S^1)\times C^\infty(S^1)\to\C
$$
given by
\begin{equation}
(f_1, f_2)\mapsto
\frac 1{2\pi}
 \int_0^{2\pi}f_1(e^{i\phi})f_2(e^{i\phi})\,d\phi
=\frac 1{2\pi}
 \int_0^{2\pi}f_1(z)f_2(z)\,\frac{dz}{z}
.
\label{eq:bilinear-pairing}
\end{equation}

\begin{observation}
\label{obs:u11-dual}
Representations $T_{p|q}$ and $T_{1-p|1-q}$
are dual with respect to $\Pi$, i.e.,
\begin{equation}
\Pi\left(T_{p|q}(g)f_1,T_{1-p|1-q}(g)f_2\right)
=\Pi(f_1,f_2)
.
\label{eq:Pi-Pi}
\end{equation}
\end{observation}

{\sc Proof.} After  simple cancelations
we get the following expression
in the left hand side of (\ref{eq:Pi-Pi})
$$
\frac 1{2\pi i}\int_{|z|=1}
f_1 \left(\frac{b+\ov a z}{a+\ov b z}\right)
f_2 \left(\frac{b+\ov a z}{a+\ov b z}\right)
\cdot(a+\ov b z)^{-1} \ov{(a+\ov b z)}^{\,\,-1}
\frac{dz}z
.
$$
Keeping in  mind $\ov z=z^{-1}$, we transform
$$
(a+\ov b z)^{-1} \ov{(a+\ov b z)}^{\,-1}
\frac{dz}z
=(a+\ov b z)^{-1} (\ov b+ az)^{\,-1}\,dz
= \left(\frac{b+\ov a z}{a+\ov b z}\right)^{-1}
d\left(\frac{b+\ov a z}{a+\ov b z}\right)
.
$$
Now the integral comes to the desired form:
$$
\qquad\qquad\qquad\qquad\qquad\qquad
\frac 1{2\pi i}\int_{|u|=1} f_1(u)\, f_2(u)\,
\frac{du} u
.
\qquad\qquad\qquad\qquad\qquad \square
$$

\SS

We also define a sesquilinear map
$$
\Pi^*:C^\infty(S^1)\times C^\infty(S^1)\to\C
$$
 by
\begin{equation}
\Pi^*(f_1, f_2):=\Pi (f_1,\ov f_2)=
 \int_0^{2\pi}f_1(z)\ov{f_2(z)}\,\frac{dz}z
.
\label{eq:second-pairing}
\end{equation}

\begin{observation}
\label{obs:u11-dual-2}
Representations $T_{p|q}$
and $T_{1-\ov q|1-\ov p}$
are dual with respect to $\Pi^*$.
\end{observation}

{\sc Proof} is same.
\hfill $\square$


\SS

{\bf\punct Intertwining operators.}
Consider  the integral operator
$$
I_{p|q}:C^\infty(S^1)\to C^\infty(S^1)
$$
 given by
\begin{equation}
\label{eq:Ipq}
I_{p|q} f(u)=\frac 1{2\pi i\,\Gamma(p+q-1)}\int_{|z|=1}
(1-z\ov u)^{\{p-1|q-1\}}
f(z)\,\frac{dz}z
,
\end{equation}
where the function $(1-z\ov u)^{\{p-1|q-1\}}$
is defined by
\begin{equation}
(1-z\ov u)^{\{p-1|q-1\}}:=
\lim_{t\to 1^-} (1-tz\ov u)^{\{p-1|q-1\}}
\label{eq:tzu}
\end{equation}
The integral converges if $\Re (p+q)>-1$.

\begin{theorem}
\label{l:pq-analytic}
The map $(p|q)\mapsto I_{p|q}$ admits the
analytic continuation to
 a holomorphic operator-valued function on $\C^2$.
\end{theorem}

\begin{theorem}
\label{th:u11-intertwine}
The operator $I_{p|q}$ intertwines
$T_{p|q}$ and $T_{1-q|1-p}$, i.e.,
$$
T_{1-p|1-q}(g)\,I_{p|q}=I_{p|q}\, T_{p|q}(g)
.
$$
\end{theorem}

\begin{corollary}
\label{cor:pqz}
If $p\notin \Z$, $q\notin \Z$,
then the representations $T_{p|q}$
and $T_{1-q|1-p}$ are equivalent.
\end{corollary}


{\bf\punct Proof of Theorems \ref{l:pq-analytic},
\ref{th:u11-intertwine}.}

\begin{lemma}
\label{l:bilat}
The expansion of the distribution (\ref{eq:tzu})
 into the Fourier
series is given by
\begin{align}
(1-z\ov u)^{p-1}(1-\ov z u)^{q-1}
=
\frac{\Gamma(p+q-1)}{\Gamma(p)\Gamma(q)}
\sum_{n=-\infty}^\infty
 \frac{(1-q)_{n}}{(p)_n}
\left(\frac zu\right)^n
\label{eq:bilat0}
=\\=
\Gamma(p+q-1)
\sum_{n=-\infty}^\infty
\frac{(-1)^n}{\Gamma(p+n)\Gamma(q-n)}
\left(\frac zu\right)^n
\label{eq:bilat}
.
\end{align}
\end{lemma}

{\sc Proof.} Let $\Re p$, $\Re q$ be sufficiently
large. Then we  write
\begin{equation}
(1-z\ov u)^{p-1}(1-\ov z u)^{q-1}=
\Bigl[ \sum_{j\ge 0} \frac{(1-p)_j}{j!}
 \left(\frac zu\right)^j\Bigr]
 \cdot
\Bigl[ \sum_{l\ge 0} \frac{(1-q)_l}{l!}
 \left(\frac uz\right)^l\Bigr]
\label{eq:1-zu}
\end{equation}
and open brackets in (\ref{eq:1-zu}).
For instance, the coefficient at $(z/u)^0$
is
$$
\sum_{k\ge0} \frac{(1-p)_k(1-q)_k}{k!\,k!}
=\,{}_2F_1(1-p,1-q;1;1)
,
$$
where $_2F_1$ is the Gauss hypergeometric function.
We  evaluate the sum with the Gauss summation
formula for $_2F_1(1)$, see \cite{HTF}, (2.1.14).
\hfill $\square$

\SS

{\sc Proof of Theorem \ref{l:pq-analytic}.}
Denote by
$$c_n:=
\frac {(-1)^n}{\Gamma(p+n)\Gamma(q-n)}
$$
 the Fourier
 coefficients in (\ref{eq:bilat}).
 Evidently, $c_n$ admits holomorphic
 continuation to the whole plane%
 \footnote{The Gamma function $\Gamma(z)$
 has simple poles at $z=0$, $-1$, $-2$, \dots and does
 not have zeros. Therefore $1/(\Gamma(p+n)\Gamma(q-n))$
 has zeros at $p=-n$, $-n-1$, \dots and at
 $q=n$, $n-1$, \dots.

 In particular, if both $p$, $q$ are integer and $q<p$,
 when $I_{p|q}=0$.}
  $\C^2$.

By \cite{HTF}, (1.18.4),
$$
\frac{\Gamma(n+a)}{\Gamma(n+b)}
\sim |n|^{a-b}\qquad\text{as $n\to\pm\infty$}
.
$$
Keeping in  mind (\ref{eq:bilat0}), we get
\begin{equation}
c_n\sim \mathrm{const}\cdot |n|^{1-p-q}
 \qquad\text{as $n\to\pm\infty$}
\label{eq:asss}
.\end{equation}
 Then
$$
I_{p|q}: z^n\mapsto c_{-n} z^n
$$
and
$$
I_{p|q}: \sum a_n z^n\mapsto \sum a_n c_{-n} z^n
.$$
Obviously, this map sends smooth functions to
smooth functions.
\hfill $\square$

\SS


{\sc Proof of Corollary \ref{cor:pqz}.} In this case,
all $c_n\ne 0$.
\hfill $\square$

{\sc Proof of Theorem \ref{th:u11-intertwine}.}
The calculation is straightforward,
\begin{multline*}
T_{1-q|1-p}(g)I_{p|q} f(u)
=\\=
\frac1{2\pi i}
(a+\ov b u)^{\{q-1|p-1\}}
\int_{|u|=1}
\left(1-\left(\frac{b+\ov a u}{a+\ov b u}\right) \ov z\right)
^{\{q-1|p-1\}} f(z)\,\frac{dz}z
.
\end{multline*}
Next, we observe
$$
(a+\ov b u)
\left(1-\left(\frac{b+\ov a u}{a+\ov b u}\right)\ov z\right)
=
(a-b\ov z)\left(1-u\left(
\frac{-\ov b+\ov a\ov z} {a-b\ov z}
\right)\right)
$$
and come to
$$
\frac1{2\pi i}
\int_{|z|=1}
\left(1-u\left(
\frac{-\ov b+\ov a\ov z} {a-b\ov z}
\right)\right)^{\{q-1|p-1\}} (a-b\ov z)^{\{q-1|p-1\}}
f(z)\,\frac{dz}z
.
$$
Now we change a variable again
$$
z=\frac{b+\ov a w}{a+\ov b w},\qquad
\ov w=\frac{-\ov b+\ov a \, \ov z}{a- b \ov z}
$$
and come to the desired expression
$$
\frac 1{2\pi i}
\int_{|w|=1}(1-u\ov w)^{\{p-1|q-1\}}
f\left(\frac{b+\ov a w}{a+\ov b w}\right)
(a+\ov b w)^{\{-p|-q\}}\frac{dw}w
.$$


{\bf\punct The unitary principal series.}

\begin{observation}
A representation $T_{p|q}$ is unitary in $L^2(S^1)$
iff
\begin{equation}
\Im p=\Im q,\qquad \Re p+\Re q=1
\label{eq:u11-condition-principal}
.\end{equation}
\end{observation}

{\sc Proof} is straightforward,
also this follows from Observation
\ref{obs:u11-dual-2}.
\hfill $\square$


\SS

\begin{figure}
$\epsfbox{sobolev.7}$

\caption{The unitary principal series in coordinates
\newline
\qquad
$h=(p-q+1)/2$, $s=\frac 1i(p+q-1)/2$.
\newline
Equivalently,
\newline
 $p=h+is$, $q=1-h+is$.
\newline
The shift $h\mapsto h+1$ does not change a representation.
Also the symmetry $s\mapsto -s$ sends
 a representation to an equivalent one. Therefore
representations of the principal series are enumerated
by  the a semi-strip $0\le h<1$, $s\ge 0$. It is more reasonable
 to think that  representations
 of the unitary principal series are enumerated by points
 of a semi-cylinder $(s,h)$, where $s\ge 0$ and $h$
is defined modulo equivalence
$h\sim h+k$, where $h\in\Z$.\label{fig:principal}}
\end{figure}


\bigskip

{\bf\large C. The complementary series}


\BS

{\bf\punct The complementary series.}
Now let
\begin{equation}
0<p<1\qquad 0<q<1
\label{eq:u11-condition-complementary}
.
\end{equation}
Consider the Hermitian form on $C^\infty(S^1)$
given by
\begin{equation}
\la f_1,f_2\ra_{p|q}=
\frac 1{(2\pi i)^2\,\Gamma(p+q-1)}
\int\limits_{|z|=1}
\int\limits_{|u|=1}
(1-z\ov u)^{\{p-1|q-1\}}
 f_1(z)\ov{f_2(u)}\,
\frac {dz} z
\frac {du} u
\label{eq:inner-complementary}
.\end{equation}
By (\ref{eq:bilat}),
\begin{equation}
\la z^n, z^m\ra_{p|q}=
\frac{1}
{\Gamma(p)\Gamma(q)}
 \frac{(1-q)_{n}}{(p)_n}
 \cdot \delta_{m,n}
\label{eq:inner-complementary-2}
.
\end{equation}

\begin{theorem}
\label{th:comp-positivity}
If
$0<p<1$, $0<q<1$, then the inner product
(\ref{eq:inner-complementary}) is positive definite.
\end{theorem}

{\sc Proof.} Indeed, in this case
 all  coefficients
$$
 \frac{(1-q)_{n}}{(p)_n}=
  \frac{(1-p)_{-n}}{(q)_{-n}}
$$
in (\ref{eq:inner-complementary})
are positive.
\hfill $\square$

\begin{theorem}
\label{th:comp}
Let $0<p<1$, $0<q<1$. Then the representation
$T_{p|q}$ is unitary with respect to the
inner product $\la\cdot,\cdot\ra_{p|q}$, i.e.,
$$
\la T_{p|q}(g) f_1, T_{p|q}(g) f_2\ra_{p|q}
=\la  f_1,  f_2\ra_{p|q}
.$$
\end{theorem}

{\sc Proof.} This follows from
Theorem \ref{th:u11-intertwine}
and Observation \ref{obs:u11-dual-2}.
Indeed,
$$
\la f_1, f_2\ra_{p|q}=\Pi^*(I_{p|q} f_1,f_2)
$$
and
\begin{multline*}
\Pi^*(I_{p|q} T_{p|q}(g) f_1, T_{p|q}(g)f_2)
=
\Pi^*( T_{1-q|1-p}(g) I_{p|q} f_1, T_{p|q}(g)f_2)
=\\=
\Pi^*(  I_{p|q} f_1, f_2)=
\la f_1, f_2\ra_{p|q}
.\end{multline*}

Keeping in mind our future purposes,
we propose  another (homotopic) proof.
Substitute
$$
z=\frac{b+\ov a z'}{a+\ov b z'},
\qquad
u=\frac{b+\ov a u'}{a+\ov b u'}
$$
to  the integral in (\ref{eq:inner-complementary}).
Applying  the identity
$$
1-\left(\frac{b+\ov a z'}{a+\ov b z'}\right)
   \ov{\left(\frac{b+\ov a u'}{a+\ov b u'}\right)}
   = (a+\ov b z')^{-1}
    (1-z'\ov u') \ov{(a+\ov b u')}^{\,-1}
,$$
we get
$$
\qquad\qquad\qquad\qquad\qquad\quad
\la T_{p|q}(g)f_1, T_{p|q}(g)f_2\ra_{p|q}
.
\qquad\qquad\qquad\qquad\qquad\quad
\square
$$

\SS

\begin{figure}
$$\epsfbox{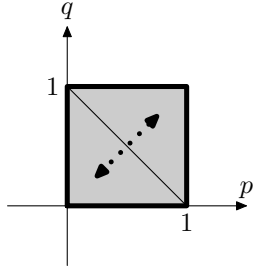}$$
\caption{The complementary series. The diagonal is contained
in the principal series (the  segment
of
the axis $Oh$ on Fig.\ref{fig:principal}). The symmetry
with respect to the diagonal sends a representation
to an equivalent representation.\label{fig:complementary}}
\end{figure}


{\bf\punct Sobolev spaces.
\label{ss:sobolev-1}}
Denote by $\cH_{p|q}$
the completion of $C^\infty(S^1)$ with respect
to  the
inner product of the complementary series.

 First, we observe that
the principal series and the complementary series
have an intersection,
see (\ref{eq:u11-condition-principal}),
(\ref{eq:u11-condition-complementary}), namely
the interval
$$
 p+q=1, \qquad 0<p<1
.$$
In this case the inner product
(\ref{eq:inner-complementary-2}) is the
$L^2$-inner product, i.e., $\cH_{p|1-p}\simeq L^2(S^1)$.

Next consider  arbitrary $(p,q)$, where $0<p<1$, $0<q<1$.
 By
(\ref{eq:asss}), the space $\cH_{p|q}$
consists of Fourier series
$\sum a_n z^n$ such that
$$
\sum_{n=-\infty}^\infty |a_n|^2
n^{1-p-q}<\infty
.$$
Thus, $\cH_{p,q}$ is the Sobolev space
$W^{(1-p-q)/2}(S^1)$.


\BS

{\bf\large D. Holomorphic and anti-holomorphic representations}

\BS


Denote by $D$ the disk $|z|<1$ in $\C$.

\SS

{\bf\punct Holomorphic (highest weight)
representations.}
Set $q=0$,
$$
T_{p|0}f(z)= f\left(\frac{b+\ov a z}{a+\ov b z}\right)
(a+\ov b z)^{-p}
.
$$
Since $|a|>|b|$, the factor $(a+\ov b z)^{-p}$
is holomorphic in the disk $D$.
Therefore the space of holomorphic functions
in $D$ is $\SU(1,1)^\sim$-invariant.
Denote the representation of  $\SU(1,1)^\sim$
 in the space of holomorphic
functions by $T^+_p$.

\begin{theorem}
a) For $p>0$ the representation $T^+_p$ is unitary,
the invariant inner product in
the space of holomorphic functions
is
\begin{equation}
\bigl\la \sum_{n\ge 0}a_nz^n,
\sum_{n\ge 0}b_nz^n\bigr\ra=
\sum_{n>0} \frac{n!}{(p)_n} a_n\ov b_n
.
\label{eq:highest-inner}
\end{equation}

b) For $p>1$ the invariant inner product
admits the following integral representation:
$$
\la f_1,f_2\ra=\frac{p-1}\pi \iint_{|z|<1}
f_1(z)\ov{f_2(z)}\,(1-|z|^2)^{p-2} d\lambda(z)
,$$
where $d\lambda(z)$ is the
Lebesgue measure in the disk.

\SS

c) For $p=1$ the invariant inner product
is
\begin{equation}
\la f_1,f_2\ra=\frac{1}{2\pi}
\int_{0}^{2\pi} f_1(e^{i\phi})\ov{f_2(e^{i\phi})}
\,d\phi
=\frac 1{2\pi i}\int_{|z|=1}
f_1(z)\ov{f_2(z)}\,\frac{dz}z
.
\label{eq:hardy-inner}
\end{equation}

\end{theorem}

We denote this Hilbert space of holomorphic
functions by $\cH_p^+$.

\SS

{\sc Proof.} The invariance of inner products
in b), c) can be easily verified by  straightforward
calculations.

To prove a), we note that weight vectors
$z^n$ must be pairwise orthogonal.

Next, operators of the Lie algebra $\frs\fru(1,1)$ must
be skew-self-adjoint.
The generators of the Lie algebra must satisfy
$$
(L_+)^*=L_-.
$$
Therefore,
$$
\la L_+ z^n,z^{n+1}\ra=\la z^{n}, L_- z^{n+1}\ra
$$
or
$$
(n+p)\la z^{n+1},z^{n+1}\ra= (n+1)\la z^n,z^n\ra
.$$
This implies a).

If $p=1$, then $\la z^n,z^n\ra=1$ for $n\ge 0$,
i.e., we get the $L^2$-inner product.
\hfill $\square$

\SS

The theorem does not provide us
 an explicit integral formula
for inner product in $\cH_p^+$ if $0<p<1$.
There is
 another
way of description of inner products in spaces
of holomorphic functions.

\SS


{\bf\punct Reproducing
kernels.}

\begin{theorem}
For each $p>0$, for any $f\in\cH_p^+$,
and for each $a\in D$
\begin{equation}
\la f(z), (1-z\ov a)^{-p}\ra=f(a)
\qquad
\text{\it (the reproducing property)}.
\label{eq:reproducing-su11}
\end{equation}
\end{theorem}

{\sc Proof.}
Indeed,
$$
\la \sum a_n z^n, \sum \frac{(p)_n}{n!} z^n \ov u^n\ra
=\sum a_n \frac{(p)_n}{n!} u^n \la z^n,z^n\ra
=\sum a_n u^n=f(u).
\qquad\square
$$

In fact, the identity (\ref{eq:reproducing-su11})
is an all-sufficient definition of the inner
product. We will not discuss this
(see \cite{FK}, \cite{Ner-notes}), and prefer
 another way.

\SS


{\bf\punct Realizations
of holomorphic representations in quotient spaces.}
Consider the representation $T_{-1|-1-p}$
of the principal series,
$$
T_{-1|-1-q} f(z)=
f\left(\frac{b+\ov a z}{a+\ov b z}\right)
(a+\ov b z)^{-1}\ov{(a+\ov b z)}^{\,-1-p}
.
$$
The corresponding
 invariant Hermitian form in $C^\infty(S^1)$,
is
\begin{equation}
\la f_1,f_2\ra_{-1|-1-p}
=\frac 1{(2\pi i)^2}
\int_{|z|=1}\int_{|u|=1} (1-\ov z u)^{-p}
f_1(z)\, \ov {f_2(u)}\,\frac{dz}z\,\frac{du}u
\label{eq:()}
.\end{equation}
(we write another pre-integral factor
in comparison with (\ref{eq:inner-complementary})).
The integral  diverges for $p>1$.
However, we can  define the inner product by
$$
\la z^n,z^n\ra=
\begin{cases}
\frac{(p)_n}{n!} \qquad &\text{if $n\ge 0$}
\\
0 \qquad & \text{if $n<0$}
\end{cases}
,$$
 the latter definition is valid for all $p>0$.

\SS

We denote by $L\subset C^\infty(S^1)$
the subspace consisting
 of series $\sum_{n<0} a_n z^n$.
This subspace is $\SU(1,1)$-invariant
and our form is nondegenerate and positive
definite on the quotient space
$C^\infty(S^1)/L$.

\SS

Next, we consider the intertwining operator
$$
\wt I_{-1|-1-p}:C^\infty(S^1) \to C^\infty(S^1)
$$
as above (but we change a normalization of the integral),
$$
\wt I_{-1|-1-p} f(u)=\frac{1}{2\pi i}
\int_{|z|=1}(1-\ov z u)^{-p} f(z)\,\frac{dz}z
$$
The kernel of the operator is $L$ and the image
consists of holomorphic functions.

\begin{observation}
a) The operator $I_{-1|-1-p}$ is a unitary operator
$$
C^\infty(S^1)/L \to \cH^+_p
.$$

b) The representation $T_{-1|-1-p}$ in $C^\infty(S^1)/L$
is equivalent to the highest weight representation
$T_p^+$
\end{observation}


{\bf\punct Lowest weight representations.}
Now set $p=0$, $q>0$.
Then operators $T_{0|q}$ preserve
the subspace consisting of
'antiholomorphic' functions
$\sum_{n\le 0} a_n z^n$. Denote
by $T^-_q$ the corresponding representation
in the space of antiholomorphic functions.
These representations are unitary.

We omit further discussion because these
representations are twins of highest weight
representations.

\SS


\begin{figure}
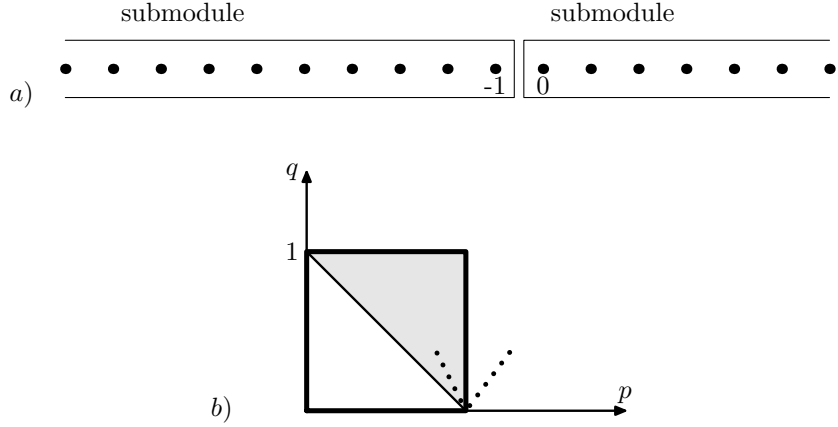


$$a)\quad \epsfbox{sobolev.5}$$

$$b)\qquad\epsfbox{sobolev.6}$$

\caption{a) The structure of the representation $T_{1|0}$.
\newline
b) Ways to $(p,q)=(1,0)$ from different
directions give origins
to different invariant Hermitian forms on $T_{1|0}$.
 By our normalization,
 the inner product is positive
definite in the gray triangle
and negative definite in the white
triangle. Therefore coming to $(1,0)$ from  the grey
triangle we get a positive form.}

\end{figure}


\BS

{\bf E. The blow-up trick}

\BS

Here we discuss a trick that produces
'unipotent' representations
of $\U(n,n)$ for $n\ge 2$,
see Subsection \ref{ss:blow-2}.

\SS


{\bf\punct The exotic case $p=1$, $q=0$.
\label{ss:blow-1}}%
In this case,
$$
T_{1|0}=T^+_1\oplus T^-_1
.
$$
Let us discuss the behavior of the
 inner
product of the complementary series
 near the point $(p|q)=(1|0)$,
\begin{equation}
\la f_1,f_2\ra_{p|q}=
\frac 1{(2\pi i)^2}\int_{|z|=1}
(1-z\ov u)^{\{p-1|q-1\}}
 f_1(z)\ov{f_2(z)}\,\,\frac {dz} z
\label{eq:inner-complementary-bis}
.\end{equation}

\SS

Consider the limit of this
expression as $p\to 1$, $q\to 0$.
The Fourier coefficients of the kernel are
the following meromorphic functions
 $$
c_n(p,q)=\frac{(-1)^n\Gamma(p+q-1)}
{\Gamma(q-n)\Gamma(p+n)}
. $$
Note that

\SS

1. $c_n(p,q)$ has a  pole
at the line $p+q=1$;

\SS

2. for $n\ge 0$, the function $c_n(p,q)$
has a zero on the line $q=0$;

\SS

3. for $n<0$, the function $c_n(p,q)$ has a  zero
at the line $p=0$.

\SS

Thus our point $(p,q)=(1,0)$ lies on the intersection
of a pole and of a zero of the function $c_n(p,q)$.
Let us substitute
$$
p=1+\epsilon s\qquad q=\epsilon t,
\qquad \text{where $s+t\ne 0$}
$$
to $c_n(p,q)$
and pass to the limit as $\epsilon\to 0$.
Recall that
\begin{equation}
\Gamma(z)=\frac{(-1)^n}{n!\,(z+n)}+O(1),\qquad
\text{as $z\to -n$, where $n=0$, $1$, $2$,\dots}
\label{eq:gamma-pole}
\end{equation}

Therefore we get
$$
\lim_{\epsilon\to 0}
c_n(1+\epsilon s,\epsilon t)
=
\begin{cases}
\frac t{t+s}\qquad & \text{if $n\ge 0$}
\\
-\frac s{t+s}\qquad & \text{if $n< 0$}.
\end{cases}
$$

In particular, for $s=0$ we get $T^+_1$-inner
product, and for $t=0$ we get $T^-_1$-inner product.
Generally,
$$
\lim_{\epsilon\to 0}\la
\sum_{n=-\infty}^\infty a_n z^n,
\sum_{n=-\infty}^\infty b_n z^n  \ra_{1+\epsilon s|\epsilon t}
=
\frac t{t+s}
\sum_{n=0}^\infty a_n\ov b_n
- \frac s{t+s}
\sum_{n=-\infty}^{-1} a_n\ov b_n
.$$
Therefore we get a one-parametric family
of invariant inner products for $T_{1|0}$.
However, all of them are linear combinations
of two basis inner products mentioned above
($t=0$ and $s=0$).


\section{Stein--Sahi representations}

\COUNTERS

Here we extend constructions of the previous section to the groups
$G:=\U(n,n)$. The analogy of the circle $S^1$ is the space
$\U(n)$ of unitary matrices.

\bigskip

{\bf\large A. Construction of representations}

\bigskip


{\bf \punct Distributions $\ell_{\sigma|\tau}$.%
\label{ss:expression}}
Let $z$ be an $n\times n$ matrix
with norm $< 1$.
For $\sigma\in\C$,
we define the function $\det(1-z)^\sigma$ by
$$
\det(1-z)^\sigma:=\det\Bigl[ 1-\sigma z+
\frac{\sigma(\sigma-1)}{2!} z^2-
  \frac{\sigma(\sigma-1)(\sigma-2)}{3!} z^3+\dots \Bigr]
.$$
Extend this function to matrices $z$ satisfying $\|z\|\le 1$,
 $\det(1-z)\ne 0$
by
$$
\det(1-z)^\sigma:=\lim_{u\to z, \,\|u\|<1}\det(1-u)^\sigma
.$$
The expression $\det(1-z)^\sigma$ is continuous
in the domain $\|z\|\le 1$ except the surface
$\det(1-z)=0$.

Denote by $\det(1-z)^{\{\sigma|\tau\}}$ the function
$$
\det(1-z)^{\{\sigma|\tau\}}:=
  \det(1-z)^\sigma \det(1-\ov z)^\tau
.$$
We define the function
$\ell_{\sigma|\tau}(g)$
on the unitary group $\U(n)$ by
\begin{equation}
\ell_{\sigma|\tau}(z): =
2^{-(\sigma+\tau)n}
\det(1-z)^{\{\sigma|\tau\}}
.
\label{eq:ell}
\end{equation}
Obviously,
\begin{equation}
\ell_{\sigma|\tau} (h^{-1}zh)=\ell_{\sigma|\tau}  (z)
\qquad\text{for $z$, $h\in\U(n)$}
\label{eq:(2.1)}
.\end{equation}

 \begin{lemma}
  Let $e^{i\psi_1}$, \dots, $e^{i\psi_n}$,
where $0\le \psi_k<2\pi$,
be the eigenvalues of $z\in\U(n)$.
Then
\begin{equation}
\ell_{\sigma|\tau}(z)=
\exp\Bigl\{\frac i2(\sigma-\tau)\sum_k(\psi_k-\pi)\Bigr\}
\prod_{k=1}^n \sin^{\sigma+\tau} \frac {\psi_k}2
.
\label{eq:ell-expression}
\end{equation}
\end{lemma}

{\sc Proof.} It  suffices to verify the statement
for diagonal matrices; equivalently we must check
the identity
$$(1-e^{i\psi})^{\{\sigma|\tau\}}=
\exp\Bigl\{\frac i2(\sigma-\tau)(\psi-\pi)\Bigr\}
\sin^{\sigma+\tau} \frac {\psi}2
.$$
We have
$$
\frac12(1-e^{i\psi})=
\exp\Bigl\{\frac i2(\psi-\pi)\Bigr\}\sin\frac \psi 2
.$$
Further, both the sides of the equality
$$
2^{-\sigma}(1-e^{i\psi})^\sigma=
\exp\Bigl\{\frac i2\sigma(\psi-\pi)\Bigr\}\sin^\sigma\frac \psi 2
,$$
 are real-analytic on $(0,2\pi)$
and the substitution $\psi=\pi$ gives 1 in  both
the sides.
\hfill $\square$

\smallskip



{\bf\punct Positivity.%
\label{ss:positivity}}
Let $\Re(\sigma+\tau)<1$.
Consider the sesquilinear form on $C^\infty(\U(n))$
given by
\begin{equation}
\langle f_1, f_2\rangle_{\sigma|\tau}=
\iint_{\U(n)\times\U(n)} \ell_{\sigma|\tau}(zu^{-1})
 f_1(z) \overline {f_2(u)}
\,d\mu(z)\,d\mu(u)
\label{eq:our-form}
.
\end{equation}

For $\sigma$, $\tau\in\R$ this form is Hermitian,
i.e.,
$$
\langle f_2, f_1\rangle_{\sigma|\tau}=
\ov{\langle f_1, f_2\rangle}_{\sigma|\tau}
$$

\begin{observation}
\label{o:meromorphic}
 For fixed $f_1$, $f_2\in C^\infty(\U(n))$,
 this expression admits
a meromorphic continuation
in $\sigma$, $\tau$ to the whole $\C^2$.
\end{observation}

This follows from general facts about distributions;
however, this  fact is a corollary of the  expansion of the
distributions $\ell_{\sigma|\tau}$ in characters,
see Theorem \ref{th:main}.
This expansion implies also the following theorem:

\begin{theorem}
\label{th:our-positivity}
 For $\sigma, \tau\in\R\setminus \Z$,
the inner product (\ref{eq:our-form})
 is positive definite (up to a sign)
iff integer parts of $-\sigma-n$ and $\tau$ are equal.
\end{theorem}

The domain of positivity is the union
of the dotted squares on Figure \ref{fig:big}.

For $\sigma$, $\tau$ satisfying this theorem, denote by
$\cH_{\sigma|\tau}$ the completion of $C^\infty(\U(n))$
with respect to our inner product.

\begin{figure}
\kartinka
\caption{Unitarizability conditions
for $\U(n,n)$. The case $n=5$. \label{fig:big}}
\end{figure}


\SS

{\bf\punct The group $\U(n,n)$.%
\label{ss:unn}}
Consider the linear space $\C^n\oplus\C^n$ equipped
with the indefinite Hermitian form
\begin{equation}
\{v\oplus w, v'\oplus w'\}=\langle v,v'\rangle_{\C^n\oplus 0}-
                             \langle w,w'\rangle_{0\oplus\C^n}
,
\label{eq:herm-cn-cn}
\end{equation}
where $\langle \cdot,\cdot\rangle$ is the standard inner product
in $\C^n$.
Denote by $\U(n,n)$ the group  of linear operators
in $\C^n\oplus\C^n$ preserving the  form $\{\cdot,\cdot\}$.
We write elements of this group as block
$(n+n)\times(n+n)$ matrices
$g:=\begin{pmatrix}
 a&b\\ c& d\end{pmatrix}$.
By definition, such matrices satisfy the condition
\begin{equation}
g\begin{pmatrix} 1&0\\0&-1\end{pmatrix}  g^*=
\begin{pmatrix} 1&0\\0&-1\end{pmatrix}
.
\label{eq:unn}
\end{equation}

\begin{lemma}
\label{l:unn-un}
The following formula
\begin{equation}
z\mapsto z^{[g]}:=(a+zc)^{-1}(b+zd)
,\qquad z\in\U(n),\,\,g=
\begin{pmatrix} a&b\\ c& d\end{pmatrix}
\in\U(n,n)
\label{eq:unn-un}
\end{equation}
determines an action of the group $\U(n,n)$  on
the {\sf space} $\U(n)$.
\end{lemma}


Proof is given in Subsection \ref{ss:proof-unn-un}.

\SS


{\bf\punct Representations $\rho_{\sigma|\tau}$
of $\U(n,n)$.%
\label{ss:rho}} Denote by $\U(n,n)^\sim$ the universal
covering of the group $\U(n,n)$, see for details
Subsection \ref{ss:universal-covering-unn}.
Fix $\sigma$, $\tau\in\C$. We define
an action of
 $\U(n,n)^\sim$   in the space $C^\infty(\U(n))$
by the linear operators
\begin{equation}
\rho_{\sigma|\tau}(g) f(z)=
f(z^{[g]})\det\nolimits^{\{-n-\tau|-n-\sigma\}}(a+zc)
.\end{equation}

 We must explain the meaning
 of the complex power in this formula. First,
$$
a+zc=(1+zca^{-1})a
$$
The defining equation (\ref{eq:unn})
implies
$\|ca^{-1}\|<1$. Hence, for all matrices $z$ satisfying
$\|z\|\le 1$, complex powers of $1+zca^{-1}$
are well defined.
Next,
$$
\det(a)^{-n-\tau|-n-\sigma}:=
\exp\Bigl\{-(n+\tau)\ln\det a-(n+\sigma)\ov{\ln\det a}
\Bigr\}
$$
It is a well-defined function on $\U(n,n)^\sim$.
We set
$$
\det(a+zc)^{-n-\tau|-n-\sigma}
:=
\det\Bigl[ (1+zca^{-1})^{-n-\tau|-n-\sigma}\Bigr]
\det(a)^{-n-\tau|-n-\sigma}
$$



{\bf\punct The Stein--Sahi representations.
\label{ss:ss}}

\begin{proposition}
\label{pr:invariance-kernel}
The operators
$\rho_{\sigma|\tau}(g)$ preserve the form
$\langle\cdot,\cdot\rangle_{\sigma|\tau}$.
\end{proposition}

Proof is given in Subsection \ref{ss:proof-invariance}.

\begin{corollary}
\label{cor:uniary-rep}
 For $\sigma$, $\tau$
satisfying the positivity conditions
of Theorem \ref{th:our-positivity},
the representation $\rho_{\sigma|\tau}$ is unitary
in the Hilbert space $\cH_{\sigma|\tau}$.
\end{corollary}


{\bf\punct The degenerate principal series.%
\label{ss:degenerate-unn}}

\begin{proposition}
\label{pr:degenerate-principal}
Let $\Re(\rho+\sigma)=-n$, $\Im \sigma=\Im \tau$.
Then the representation $\rho_{\sigma|\tau}$
is unitary in $L^2(\U(n))$.
\end{proposition}

\begin{proposition}
\label{pr:L2}
$$
\left.
\frac
{\la f_1,f_2\ra_{\sigma|\tau} }
{\prod\nolimits_{j=1}^n \Gamma(\sigma+\tau+j)}
\right|_{\sigma=-n-\tau}
=\const\cdot \int_{\U(n)} f_1(u)\ov{f_2(u)}\,d\mu(u)
$$
.
\end{proposition}


 {\bf\punct Shifts of parameters.%
\label{ss:shift}}

\begin{proposition}
\label{pr:shift}
 For integer $k$,
$$
\rho_{\sigma+k|\tau-k}\simeq (\det g)^k
\cdot\rho_{\sigma|\tau}
$$
 The intertwining operator is
the multiplication by the determinant
$$
F(z)\mapsto F(z)\det(z)^k
.$$
This operator also defines an isometry
of the corresponding Hermitian forms.
\end{proposition}

\begin{figure}

\epsfbox{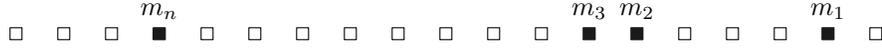}

\caption{A '{\it Maya diagram}' for signatures.
We draw the integer 'line' and fill the boxes
$m_1$, \dots, $m_n$ with black.}

\end{figure}


\bigskip

{\bf\large B. Expansions of distributions $\ell_{\sigma|\tau}$
in characters. Positivity}

\medskip


{\bf\punct Characters of $\U(n)$.%
\label{ss:characters}} See Weyl's book \cite{Wey}.
 The set of finite dimensional representations
of $\U(n)$ is parameterized by collections of integers
({\it signatures})
$$
\bm: \quad m_1> m_2 >\dots > m_n
.$$
The character $\chi_\bm$ of the  representation%
\footnote{Explicit constructions of representations
of $\U(n)$ are not used below.} $\pi_\bm$
(a {\it Schur function})
corresponding to a signature $\bm$
is given by
\begin{equation}
\chi_\bm(z)=
   \frac
    {\det_{k,j=1,2,\dots,n} \bigl\{e^{i m_j \psi_k}\bigr\} }
            {\det_{k,j=1,2,\dots,n} \bigl\{e^{i(j-1) \psi_k}\bigr\}    }
,
\label{eq:(2.3)}
\end{equation}
where $e^{i\psi_k}$ are the eigenvalues of $z$.
  Recall that the denominator admits the decomposition
\begin{equation}
 \det_{k,j} \bigl\{e^{i(j-1) \psi_k}\bigr\}
=\prod\nolimits_{l<k} (e^{i\psi_l}-e^{i\psi_k} ).
\label{eq:(2.4)}
\end{equation}

The dimension  of $\pi_m$ is
\begin{equation}
\dim \pi_\bm=\chi_\bm(1)=
\frac{\prod_{0\le\alpha<\beta\le n} (m_\alpha-m_\beta)}
{\prod_{j=1}^n j!}
.
\label{eq:(2.5)}
\end{equation}


{\bf\punct Central functions.%
\label{ss:central}}
A function $F(z)$ on $\U(n)$ is called {\it central}
if
$$
F(h^{-1} z h)=F(z) \qquad\text{for all $z$, $h\in\U(n)$}.
$$
In particular  characters and
$\ell_{\sigma|\tau}$ are central functions.

For  central functions $F$ on $\U(n)$,
 the following {\it Weyl integration formula}
 holds
\begin{multline}
\int\limits_{\U(n)} F(z) \,d\mu(z)
=
\frac 1{(2\pi)^n n!}
\int\limits_{0<\psi_1<2\pi}\dots\int\limits_{0<\psi_n<2\pi}
F\bigl(
\mathrm{diag}(e^{i\psi_1},\dots,e^{i\psi_n})
\bigr)
\times\\ \times
\Bigl| \prod_{m<k} (e^{i\psi_m}-e^{i\psi_k} ) \Bigr|^2
\,\prod_{k=1}^n d\phi_k
,
\label{eq:(2.6)}
\end{multline}
where $\mathrm{diag}(\cdot)$
is a diagonal matrix with given entries.

Any  central function $F\in L^2(\U(n))$
 admits an expansion in characters,
$$
F(z)=\sum\nolimits_\bm c_\bm\chi_m(z)
,$$
where the summation is given over all signatures $\bm$
and  the coefficients
 $c_\bm$ are  $L^2$-inner products
$$
c_\bm=\int_{\U(n)} F(z) \overline{\chi_\bm(z)}\,d\mu(z)
.$$

Note that $\ov\chi_\bfm=\chi_{\bfm^*}$,
where
$$
\bfm^*:=(n-1-m_n,\dots,n-1-m_2,n-1-m_1)
$$

Applying formula (\ref{eq:(2.6)}), explicit expression
(\ref{eq:(2.3)}) for characters,
and formula (\ref{eq:(2.4)}) for the denominator, we obtain
\begin{multline}
c_{\bm}=
\frac 1{(2\pi)^n n!}
\int\limits_{0<\psi_1<2\pi}\dots\int\limits_{0<\psi_n<2\pi}
F\Bigl(\mathrm{diag}\bigl\{e^{i\psi_1},\dots,e^{i\psi_n}\bigr\}\Bigr)
\times\\ \times
 {\det_{k,j=1,2,\dots,n} \bigl\{e^{i(j-1) \psi_k}\bigr\}    }
      \det_{k,j=1,2,\dots,n} \bigl\{e^{-i m_j \psi_k}\bigr\}
\,\prod_{k=1}^n d\phi_k
.
\label{eq:(2.7)}
\end{multline}

Let $F(z)$ be multiplicative with respect to eigenvalues,
$$
F(z)=\prod_k f\bigl(e^{i\phi_k}\bigr)
$$
(for, instance $F=\ell_{\sigma|\tau}$,
 see (\ref{eq:ell-expression})).
Then
we can apply the following
 simple lemma
 (see e.g. \cite{Ner-determinant}).

\begin{lemma}
\label{l:int-det}
 Let $X$ be a set,
\begin{multline}
\int_{X^n} \prod_{k=1}^n f(x_k) \, \det\limits_{k,l=1,\dots n}
\{u_l(x_k)\} \det\limits_{k,l=1,\dots n} \{v_l(x_k)\}
 \prod_{j=1}^n
dx_j
=\\
= n! \det\limits_{l,m=1,\dots, n}\Bigl\{ \int_X f(x)u_l(x)v_m(x)\,dx
\Bigr\}
\label{eq:determinant-lemma}
\end{multline}
.
\end{lemma}

\SS


{\bf\punct Lobachevsky beta-integrals.
\label{ss:lobachevsky}}
We wish to apply  Lemma \ref{l:int-det}
  to
 functions $\ell_{\sigma|\tau}$.
For this purpose, we need for the following
integral, see \cite{GR}, 3.631,1,  3.631,8,
\begin{equation}
\int_0^\pi \sin^{\mu-1} (\phi) \, e^{ib\phi}\,d\phi=
 \frac{2^{1-\mu}\pi
\Gamma(\mu) e^{ib\pi/2}}
    {\Gamma\bigl((\mu+b+1)/2\bigr)
     \Gamma\bigl((\mu-b+1)/2\bigr) }
\label{eq:Lob}
.
\end{equation}

It is equivalent to the identity (\ref{eq:bilat}).

In a certain sense, the integral
(\ref{eq:Kad}) is a multivariate analog
 of the Lobachevsky integral. On the other hand,
(\ref{eq:Kad}) is a special case of the  modified
Kadell integral \cite{Ner-stein}.

\SS


{\bf\punct Expansion of the function $\ell_{\sigma|\tau}$
in characters.\label{ss:expansion}}

\smallskip

\begin{theorem}
\label{th:main}
 Let $\Re (\sigma+\tau)<1$.
Then
\begin{align}
\ell_{\sigma|\tau}(g)&= \nonumber
\\
&=
\frac{(-1)^{n(n-1)/2}\sin^n(\pi\sigma)
2^{-(\sigma+\tau)n}}{\pi^n}
\prod_{j=1}^n\Gamma(\sigma+\tau+j)
\times  \notag \\
&\qquad\times
\sum\limits_\bm\Biggl\{
\prod\limits_{1\le\alpha<\beta\le n} (m_\alpha-m_\beta)
\prod\limits_{j=1}^n\frac{\Gamma(-\sigma+m_j-n+1)}
    {\Gamma(\tau+m_j+1)}\chi_\bm(g)
\Biggr\}
\label{eq:expansion-1}
=\\
&=
(-1)^{n(n-1)/2}
2^{-(\sigma+\tau)n}
\prod_{j=1}^n\Gamma(\sigma+\tau+j)
\times\notag
\\ &\quad\times
\sum\limits_\bm\Biggl\{
\frac{(-1)^{\sum m_j}
\prod\limits_{1\le\alpha<\beta\le n} (m_\alpha-m_\beta)}
{\prod\limits_{j=1}^n\Gamma(\sigma-m_j+n)
    \Gamma(\tau+m_j+1)}\chi_\bm(g) \Biggr\}
.
\label{eq:expansion-2}
\end{align}
\end{theorem}

Proof is contained in
Subsection \ref{ss:proof}.
For the calculation we need for Lemma \ref{l:krat}
proved in the next subsection.

\SS


{\bf\punct A determinant identity.
\label{ss:krat}}
Recall that the {\it  Cauchy determinant}
(see e.g. \cite{Kra})
is given by
\begin{equation}
\det\limits_{kl}\left\{ \frac 1{x_k+y_l}\right\}=
\frac{
\prod_{1\le k <l\le n}(x_k-x_l)
\cdot
\prod_{1\le k <l\le n}(y_k-y_l)}
{\prod\limits_{1\le k, l\le n} (x_k+y_l)}
.
\label{eq:cauchy}
\end{equation}

The following version of the Cauchy determinant
 is also well known.

\smallskip

\begin{lemma}
\label{l:cauchy-2}
\begin{multline}
\det\begin{pmatrix}
1&1&1& \dots & 1 \\
\frac 1{x_1+b_1} & \frac 1{x_2+b_1} & \frac 1{x_3+b_1}&
    \dots & \frac {1\vphantom{1^G}}
{x_n+b_1}\\
\frac {1\vphantom{1^G}}{x_1+b_2} & \frac 1{x_2+b_2} & \frac 1{x_3+b_2}&
    \dots & \frac 1{x_n+b_2}\\
\vdots&\vdots&\vdots&\ddots &\vdots\\
\frac {1\vphantom{1^G}}{x_1+b_{n-1}} & \frac 1{x_2+b_{n-1}} & \frac 1{x_3+b_{n-1}}&
    \dots & \frac 1{x_n+b_{n-1}}\\          \end{pmatrix}
=\\=
\frac{
\prod_{1\le k <l\le n} (x_k-x_l)
 \prod_{1\le \alpha<\beta\le n-1} (b_\alpha-b_\beta)}
{\prod\limits_{\begin{smallmatrix} 1\le k \le n\\
                   1\le \alpha \le n-1 \end{smallmatrix}}
(x_k+b_\alpha)}
.
\label{eq:cauchy-2}
\end{multline}
\end{lemma}

{\sc Proof.}  Let $\Delta$ be the Cauchy determinant
(\ref{eq:cauchy}).
Then
$$
y_1\Delta=\begin{pmatrix}
\frac {y_1}{x_1+y_1}&  \frac {y_1}{x_2+y_1}&\dots & \frac {y_1}{x_1+y_1}\\
\frac 1{x_1+y_2}&  \frac 1{x_2+y_2}&\dots & \frac 1{x_n+y_2}\\
\vdots&\vdots&\ddots &\vdots\\
\frac 1{x_1+y_n}&  \frac 1{x_2+y_n}&\dots & \frac 1{x_n+y_n}
\end{pmatrix}
.$$
We take  $\lim\limits_{y_1\to \infty} y_1 \Delta$
and substitute $y_{\alpha+1}=b_\alpha$.
\hfill $\square$

\smallskip

The following determinant is
a rephrasing of   \cite{Kra}, Lemma 3.

\begin{lemma}
\label{l:krat}
\begin{multline}
\!\!\!\!\!\!\!\!\!
\det\begin{pmatrix}
1&1&1& \dots & 1 \\
\frac {x_1+b_1\palka} {x_1+a_1} & \frac {x_2+a_1}{x_2+b_1} &
     \frac {x_3+a_1}{x_3+b_1}&
    \dots & \frac {x_n+a_1}{x_n+b_1}\\
\frac {(x_1+a_1)\palka(x_1+a_2)}{(x_1+b_1)(x_1+b_2)}
  & \frac {(x_2+a_1)(x_2+a_2)}{(x_2+b_1)(x_2+b_2)}
  & \frac {(x_3+a_1)(x_3+a_2)}{(x_3+b_1)(x_3+b_2)}&
    \dots & \frac {(x_n+a_1)(x_n+a_2)}{(x_n+b_1)(x_n+b_2)}\\
\vdots & \vdots & \vdots &\ddots &\vdots\\
\frac {\prod\limits_{1\le m\le n-1} (x_1+a_m)}
{\prod\limits_{1\le m\le n-1} (x_1+b_m)}&
\frac {\prod\limits_{1\le m\le n-1} (x_2+a_m)}
     {\prod\limits_{1\le m\le n-1} (x_2+b_m)}&
\frac {\prod\limits_{1\le m\le n-1} (x_3+a_m)}
   {\prod\limits_{1\le m\le n-1} (x_3+b_m)} &
\dots&
\frac {\prod\limits_{m:\,1\le m\le n-1} (x_n+a_m)}
 {\prod\limits_{m:\,1\le m\le n-1} (x_n+b_m)}
\end{pmatrix}
=\\=
\frac{\prod\limits_{1\le k<l\le n} (x_k-x_l)
  \prod\limits_{1\le \alpha\le\beta\le n-1} (a_\alpha-b_\beta)}
{\prod\limits_{1\le k\le n, 1\le\beta\le n-1} (x_k+b_\beta)}
.
\label{eq:krat}
\end{multline}
\end{lemma}

{\sc Proof.}
Decomposing a matrix element into a sum
of partial fractions, we obtain
\begin{multline*}
\frac{(x_k+a_1)\dots (x_k+a_\alpha)}
     {(x_k+b_1)\dots (x_k+b_\alpha)}=
1+\sum\limits_{1\le \beta \le \alpha}
\frac{\prod_{j\le\alpha}(a_j-b_\beta)}
     {  \prod_{j\le\alpha, j\ne\beta}(b_j-b_\beta)}
\cdot \frac 1{x_k+b_\beta}
\end{multline*}


Therefore the  $(\alpha+1)$-th row
is a linear combination of the following rows:
$$
\begin{matrix}
\Bigl(&1&1&\dots&1&\Bigr), \\
\Bigl(&\frac 1{x_1+b_1}& \frac 1{x_2+b_1}
 &\dots&\frac 1{x_n+b_1} &\Bigr),\\
\dots&\dots&\dots&\dots&\dots&\dots\\
\Bigl(&\frac 1{x_1+b_\alpha}&
 \frac 1{x_2+b_\alpha} &\dots&\frac 1{x_n+b_\alpha} &\Bigr).\\
\end{matrix}
$$
Thus our determinant equals
$$
\prod_{\alpha=1}^{l-1}
\frac{\prod_{j=1}^{\alpha}(a_j-b_\alpha)}
{\prod_{j=1}^{\alpha-1}(b_j-b_\alpha)}
\cdot
\det
\begin{pmatrix}
1&1&\dots&1\\
\frac 1{x_1+b_1}& \frac 1{x_2+b_1} &\dots&\frac 1{x_n+b_1} \\
\vdots&\vdots&\ddots&\vdots\\
\frac 1{x_1+b_\alpha}&
 \frac 1{x_2+b_\alpha} &\dots&\frac 1{x_n+b_\alpha}\\
\end{pmatrix}
.$$
and we refer to Lemma \ref{l:cauchy-2}.
\hfill $\square$

\SS


{\bf\punct Proof of Theorem \ref{th:main}.
\label{ss:proof}}
We must evaluate the
inner product
$$
\int_{\U(n)} \ell_{\sigma|\tau}(g)\,\ov{\chi_\bm(g)}\,d\mu(g)
.
$$
Applying (\ref{eq:(2.7)}), we get
\begin{multline}
\frac1{(2\pi)^n\,n!}
\int\limits_{0<\psi_k<2\pi}
\prod\limits_{j=1}^n
\Bigl[\sin^{\sigma+\tau}
  \bigl( \psi_j/ 2\bigr)\cdot
 \exp\Bigl\{\frac i2(\sigma-\tau)(\psi_j-\pi)\Bigr\}
\Bigr]
\times\\
\times \det\limits_{1\le k,l\le n}
    \{e^{-im_k\psi_l}\}
\cdot\det\limits_{1\le k,l\le n}
    \{e^{i(k-1)\psi_l}\}
\prod\limits_{l=1}^n d\psi_l
\label{eq:Kad}
.\end{multline}
By Lemma \ref{l:int-det},
we reduce this integral to
$$
\frac1{(2\pi)^n}
\det\limits_{1\le k,j\le n} I(k,j)
,$$
where
$$
I(k,j)=e^{-i(\sigma-\tau)\pi/2}
\int_0^{2\pi}
\sin^{\sigma-\tau}( \psi/ 2)
\cdot
 \exp\bigl\{i(\,(\sigma+\tau)/2+k-1-m_j)\bigr\}
\,d\psi
.$$
We apply the Lobachevsky integral (\ref{eq:Lob})
and get
$$
I(k,j)=\frac{2^{1-\sigma-\tau}\pi\Gamma(\sigma+\tau+1)\, (-1)^{k-1-m_j} }
            {\Gamma(\sigma+k-m_j)\Gamma(\tau-k+m_j+2)}
$$
Applying standard formulas for $\Gamma$-function,
we come to
\begin{multline*}
I(k,j)=2^{1-\sigma-\tau}\Gamma(\sigma+\tau+1)\, \sin(-\sigma\pi)
\cdot
           \frac{ \Gamma(-\sigma+m_j-k+1)}
            {\Gamma(\tau+m_j-k+2)}
                 =\\=
2^{1-\sigma-\tau}\Gamma(\sigma+\tau+1)\, \sin(-\sigma\pi)
\cdot
           \frac{ \Gamma(-\sigma+m_j-n+1)}
                  {\Gamma(\tau+m_j-n+2)}
\, \cdot
       \,  \boxed {  \frac{ (-\sigma+m_j-n+1)_{n-k}}
                  { (\tau+m_j-n+2)_{n-k}}
               }
\end{multline*}
The factors outside the box do not depend on
on $k$. Thus, we must evaluate the determinant
$$
\det\limits_{1\le k,j\le n}
           \frac{ (-\sigma+m_j-n+1)_{n-k}}
                  { (\tau+m_j-n+2)_{n-k}}
.$$
Up to a permutation of rows, it is a determinant of
the form described in Lemma \ref{l:krat}
with
$$x_j=m_j,\qquad a_j=-\sigma-n+j,\qquad b=\tau-n+j+1
.$$
After a  rearrangement of the factors,
we obtain the required result.
\hfill $\square$

\smallskip


{\bf\punct Characters of compact groups. Preliminaries.
\label{ss:compact}}
First, recall some standard facts on
characters of compact groups, for details,
see e.g. \cite{Kir}, 9.2, 11.1.

Let $K$ be a compact Lie group
equipped with the Haar measure $\mu$, let $\mu(K)=1$.
Let
$\pi_1$, $\pi_2$, \dots
be the complete collection of pairwise distinct
irreducible representations of $K$.
Let $\chi_1$, $\chi_2$, \dots be their characters.
Recall the orthogonality relations,
\begin{equation}
\langle \chi_k,\chi_l \rangle_{L^2(K)}=
\int_K \chi_k(h)\ov{\chi_l(h)}\,d\mu(h)=
\delta_{k,l}
\label{eq:ortho}
\end{equation}
and
\begin{equation}
\chi_k*\chi_l= \begin{cases}
\frac 1{\dim\pi_k} \chi_k &\qquad\text{if $k=l$},\\
                 0 &\qquad \text{if $k\ne l$}.
\end{cases}
\label{eq:convolution}
\end{equation}
where $*$ denotes the convolution on the group,
$$
u*v(g)=\int_K u(gh^{-1})\,v(h)\,d\mu(h).
$$

Consider the action of the group $K\times K$
in $L^2(K)$ by the left and right shifts
$$(k_1,k_2):\,\,f(g)\mapsto f(k_1^{-1}gk_2)
.
$$

The representation of $K\times K$
in $L^2(K)$ is a multiplicity free direct sum of
irreducible representations having the form
$\pi_k\otimes\pi_k^*$, where $\pi_k^*$
denotes the dual representation,
\begin{equation}
L^2(K)\simeq \bigoplus_k \pi_k\otimes \pi^*_k
\label{eq:frobenius}
.\end{equation}
Denote
by $V_k\subset L^2(K)$  the space of representation
$\pi_k\otimes \pi^*_k$.
Each distribution $f$ on $K$ is a sum
of 'elementary harmonics',
$$
f=\sum\nolimits_k f^{k},\qquad f_k\in V_k.
.
$$

The projector to a subspace $V_k$
is the convolution with the corresponding character,
\begin{equation}
f^{k}=\frac 1{\dim \pi_k} f*\chi_k
\label{eq:projector-harmonics}
\end{equation}
(in particular, $f^{k}$ is smooth).

\begin{observation}
\label{obs:rapid}
 Let $f$ be a function on $\U(n)$,
$f=\sum_{\bfm} a_m f^{\bfm}$, where
$
f^{\bfm}\in V_{\bfm}
.$

a) $f\in C^\infty(\U(n))$ iff
$$
\|f^{\bfm}\|_{L^2}=o\left(\sum m_j^2\right)^{-L}
\qquad\text{for all $L$}
.$$

b) $f$ is a distribution on $\U(n)$
iff there exists $L$ such that
$$
\|f^{\bfm}\|_{L^2}=o\left(\sum m_j^2\right)^{L}
.$$
\end{observation}

Proof: Note that $f\in L^2(\U(n))$ iff
$\sum\|f^\bfm\|^2_{L^2}<\infty$.
Denote by $\Delta$ be the second order invariant
 Laplace operator
on $\U(n)$. Then $\Delta f^{\bfm}=q(\bfm)f^{\bfm}$,
where $q(\bfm)=\sum m_j^2 +\dots$ is an explicit quadratic expression
 in $\bfm$. For $f\in C^\infty$ we have $\Delta^p f\in C^\infty$;
this implies the first statement.
Since $q(\bfm)$ has a finite number of zeros (one),
  the second statement follows
from a) and the duality.
 \hfill $\square$

\SS


{\bf\punct Hermitian forms defined by kernels.%
\label{ss:forms-kernels}}
 Let $\Xi$ be  a central distribution on $K$
satisfying $\Xi(g^{-1})=\ov{\Xi(g)}$.
Consider the following Hermitian form
on $C^\infty(K)$
\begin{equation}
\langle f_1,f_2\rangle=
\iint_{K\times K} \Xi(gh^{-1})\,f_1(h)\ov{f_2(g)}
\,d\mu(h)\,d\mu(g)
.
\label{eq:xi-f-f}
\end{equation}
Consider the expansion of $\Xi$ in characters
$$\Xi=\sum_k c_k\chi_k.$$

\begin{lemma}
\label{l:expansion-product}
\begin{equation}
\langle f_1, f_2\rangle =
\sum_k\frac{c_k}{\dim\pi_k}
\int_{\U(n)} f_1^{k}(h)\ov{ f_2^{k}(h) }
\,d\mu(h)
\label{eq:expansion-product}
.\end{equation}
\end{lemma}

{\sc Proof.} The Hermitian form (\ref{eq:xi-f-f})
 is
$K\times K$-invariant. Therefore the subspaces
$V_k\simeq\pi_k\otimes \pi_k^*$ must be pairwise
orthogonal. Since $\pi_k\otimes \pi_k^*$
is an irreducible representation of $K\times K$,
it admits a unique up to a factor $K\times K$-invariant
Hermitian form. Therefore
it is sufficient to find
these  factors.

 Set $f_1=f_2=\chi_k$.
We evaluate
$$
\iint_{K\times K} \left(\sum_k c_k \chi_k(gh^{-1})\right)\,
     \chi_{k}(h) \ov{\chi_k(g)}\,d\mu(g)\,d\mu(h)
=\frac{c_k}{\dim\pi_k}
$$
using (\ref{eq:ortho}) and (\ref{eq:convolution}).
\hfill $\square$

\smallskip


{\bf\punct Positivity.
\label{ss:proof-positivity}}
Let $\Re(\sigma+\tau)<1$.
Consider the sesquilinear form on $C^\infty(\U(n))$
given by
\begin{equation}
\langle f_1, f_2\rangle_{\sigma|\tau}=
\iint_{\U(n)\times\U(n)} \ell_{\sigma|\tau}(zu^{-1})
 f_1(z) \overline {f_2(u)}
\,d\mu(z)\,d\mu(u)
,
\label{eq:our-form-1}
\end{equation}
where the distribution $\ell_{\sigma|\tau}$ is the same as above.

\begin{observation}
\label{obs:meromorphic}
 For fixed $f_1$, $f_2\in C^\infty(\U(n))$,
 the expression $\langle f_1, f_2\rangle_{\sigma|\tau}$
 admits
a meromorphic continuation
in $\sigma$, $\tau$ to the whole $\C^2$.
\end{observation}

{\sc Proof.} Expanding $f_1$, $f_2$
 in elementary  harmonics
$$
f_1(z)=\sum_\bm f_2^\bm(z),\qquad f_2(z)=\sum_\bm f_2^\bm(z)
,$$
we get (see Lemma \ref{l:expansion-product})
$$
\langle f_1, f_2\rangle_{\sigma|\tau}=
\sum_\bm \frac{c_\bm}{\dim \pi_\bm}
\int_{\U(n)}f_1^\bm(z) \ov{ f_2^\bm(z)}\, d\mu(z)
,$$
where the meromorphic expressions for $c_\bm$ were obtained
in Theorem \ref{th:main}. The coefficients
$c_\bfm$ have polynomial growth in $\bfm$.
On the other hand, $\|f_j^\bm\|$ rapidly decrease,
see Observation \ref{obs:rapid}.
Therefore, the series converges.
\hfill $\square$

\SS

{\sc Proof of positivity. Corollary \ref{cor:uniary-rep}} We look at
the expression (\ref{eq:expansion-1}).
It suffices to examine the factor
\begin{equation}
\frac{\Gamma(-\sigma-n+m_j+1)}
    {\Gamma(\tau+m_j+1)}
,
\label{eq:Gamma/Gamma}
\end{equation}
because signs of all the remaining factors are independent
on $m_j$. Let $n\in\Z$ and $\alpha\in(0,1)$.
Then
$$
\mathrm{sign}\,\Gamma(n+\alpha)=
\begin{cases}
+1,\qquad\text{if $n\ge>0$,}
\\
(-1)^n,\qquad\text{if $n<0$}
\end{cases}
$$
Therefore (\ref{eq:Gamma/Gamma}) is positive
whenever integer parts of $\tau$ and $-\sigma-n$ equal.
\hfill $\square$


\SS

{\bf\punct The $L^2$-limit. Proof of Proposition
\ref{pr:L2}.} Thus, let $\sigma+\tau=-n$.
Then
$$
\Bigl( \prod_{j=1}^n
\Gamma(\sigma+\tau+j)\Bigr)^{-1}\ell_{\sigma|\tau}=
\const\cdot\sum (\dim \pi_\bfm) \chi_\bfm
$$
Indeed, in this case $\Gamma$-factors
in (\ref{eq:expansion-1}) cancel,
and we use
(\ref{eq:(2.5)}).

Keeping in mind (\ref{eq:expansion-product}),
we get Proposition \ref{pr:L2}.

\BS

{\bf\large C. Other proofs}

\BS

Here we prove that the operators $\rho_{\sigma|\tau}$
preserve the inner product determined by the distribution
$\ell_{\sigma|\tau}$.

\SS


{\bf\punct The universal covering of the group $\U(n,n)$.
\label{ss:universal-covering-unn}}
The fundamental group of $\U(n,n)$ is%
\footnote{By a general theorem,
 a real reductive Lie group $G$ admits a deformation retraction
to its maximal compact subgroup $K$. In our case,
$K=\U(n)\times\U(n)$ and $\pi_1(\U(n))=\Z$.
}
$$
\pi_1\bigl(\U(n,n)\bigr)\simeq\Z\oplus\Z
.
$$

The universal covering $\U(n,n)^\sim$ of $\U(n,n)$
can be identified with the the  set $\frU$ of triples
$$
\left\{\,\,\,
 \begin{pmatrix} a&b\\c&d\end{pmatrix},\,
s,\,t\right\}
\in \U(n,n)\times \C\times\C
$$
satisfying the conditions
$$
\det (a)=e^s,\qquad \det (d)= e^t.
$$
The multiplication of triples is given by the formula
$$
(g_1,s_1,t_1)\circ(g_2,s_2,t_2)=
\bigl(g_1g_2,s_1+s_2+c^+(g_1,g_2),t_1+t_2+c^-(g_1,g_2)\bigr)
,
$$
where the {\it Berezin cocycle} $c^\pm$ is given by
$$
c^+(g_1,g_2)=\tr\ln (a_1^{-1}a_3a_2^{-1}),\qquad
c^-(g_1,g_2)=\tr\ln (d_1^{-1}d_3d_2^{-1})
;$$
here $g_3=g_1g_2$, and
$g_j=\begin{pmatrix}a_j&b_j\\c_j&d_j\end{pmatrix}$.
It can be shown that $\|a_1^{-1}a_3a_2^{-1}-1\|<1$,
therefore the logarithm is well defined. On the other hand,
$$
e^{s_3}=e^{s_1+s_2+c^+(g_1,g_2)}=\det(a_1) \det(a_2)
\det(a_1^{-1}a_3a_2^{-1})=\det(a_3)
$$
This shows that the $\frU$ is closed with respect
to multiplication.

For details, see \cite{Ner-notes}.

In particular, $\det(a)$ is a well-defined single-valued
function on $\U(n,n)^\sim$. In our notation, it is given by
$$
(g,s,t)\mapsto s
.
$$


{\bf\punct Another model of $\U(n,n)$.
\label{ss:another-unn}}
We can  realize $\U(n,n)$ as the group of
$(n+n)\times(n+n)$-matrices
$g=\begin{pmatrix}\alpha&\beta\\ \gamma&\delta\end{pmatrix}$
satisfying the condition
\begin{equation}
g\begin{pmatrix} 0&i\\-i&0\end{pmatrix} g^*=
\begin{pmatrix} 0&i\\-i&0\end{pmatrix}
\label{eq:unn-another}
\end{equation}


{\bf\punct Action of $\U(n,n)$ on the space $\U(n)$.
Proof of Lemma \ref{l:unn-un}
\label{ss:proof-unn-un}}.
We must show that for
$$ z\in\U(n)
\quad \text{and} \quad
g=\begin{pmatrix} a&b\\c&d \end{pmatrix}
\in\U(n,n)
$$
we have
\begin{equation}
 z^{[g]}:=(a+zc)^{-1}(b+zd)\in \U(n,n)
\label{eq:unn-un-1}
\end{equation}

For $z\in\U(n)$, consider its graph
 $graph(z)\subset\C^n\oplus\C^n$.
It is an $n$-dimensional linear subspace, consisting of
all vectors $v\oplus vz$, where
a vector-row $v$ ranges in $\C^n$.
Since $z\in\U(n)$,
the subspace $graph(z)$ is isotropic%
\footnote{A subspace $V$ in a linear space is
{\it isotropic} with respect to an Hermitian
 form $Q$ if
$Q$ equals 0 on $V$.}
with respect to the Hermitian form
$\begin{pmatrix}1&0\\0&-1 \end{pmatrix}$.
Conversely, any $n$-dimensional isotropic subspace
in $\C^n\oplus\C^n$ is a graph
of a unitary   operator $z\in\U(n)$.

Thus we get a one-to-one correspondence between
the group $\U(n)$ and the Grassmannian of $n$-dimensional isotropic
subspaces in $\C^n\oplus\C^n$.

The group $\U(n,n)$ acts on the Grassmannian
and therefore $\U(n,n)$
acts on the space $\U(n)$.
Then (\ref{eq:unn-un-1})
 is the explicit expression for the
 latter action. Indeed
$$
\Bigl(v\oplus vz\Bigr)
\begin{pmatrix} a&b\\c&d \end{pmatrix}
=v(a+zc)\,\oplus\, v(b+zd)
$$
We denote $\xi:=v(a+zc)$ and get
$$
\xi\,\oplus\,\xi(a+zc)^{-1}(b+zd)
$$
and this completes the proof  of Lemma \ref{l:unn-un}.
\hfill $\square$

\SS

Thus, $\U(n)$ is a $\U(n,n)$-homogeneous space.
We describe without proof (it is a simple
exercise) the stabilizer of a point
$z=1$. It  is a maximal parabolic
subgroup.

In the model (\ref{eq:unn-another}) it can be
realized as the subgroup  of matrices having the structure
$$
\begin{pmatrix}\alpha&0\\\beta&\alpha^{*-1}\end{pmatrix}
$$
It is a semidirect product of $\GL(n,\C)$ and
the Abelian group $\R^{n^2}$.

In our basic model the stabilizer of $z=1$
is the semi-direct product of two  subgroups
\begin{equation}
\frac 12
\begin{pmatrix}
 \alpha+\alpha^{*-1}&\alpha-\alpha^{*-1}\\
                \alpha-\alpha^{*-1}&\alpha+\alpha^{*-1}
,\end{pmatrix}
\quad\text{where $g\in\GL(n,\C)$}
,
\label{eq:stabilizer-1}
\end{equation}
and
\begin{equation}
\begin{pmatrix}
1+iT&iT\\
-iT&1-iT
\end{pmatrix},\quad
\text{where $T=T^*$}
.
\label{eq:stabilizer-2}
\end{equation}


{\bf\punct The Jacobian.
\label{ss:proof-Jacobian}}

\begin{lemma}
\label{l:Jacobian}
 For the Haar measure $\mu(z)$ on
$\U(n)$, we have
\begin{equation}
\mu\bigl(z^{[g]}\bigr)=
|\det\nolimits^{-2n}(a+zc)| \cdot \mu(z)
.\label{eq:Jacobian}
\end{equation}
\end{lemma}

{\sc Proof.} A verification of this formula is straightforward,
we only outline the main steps. First,
$J(g,z):=|\det\nolimits^{-2n}(a+zc)|$ satisfies the chain
rule (\ref{eq:chain-rule}). Next, the formula
(\ref{eq:Jacobian}) is valid
for $g\in\U(n,n)$ having the form
$\begin{pmatrix}a&0\\0&d\end{pmatrix}$,
where $u$, $v\in\U(n)$. Indeed, the corresponding transformation
of $u\mapsto u^{[h]}$ is $u\mapsto a^{-1}ud$, its Jacobian
is 1.

Therefore we can set $z=1$, $z^{[g]}=1$.
Now we must evaluate
the determinants of the differentials
of maps $z\mapsto z^{[g]}$ at $z=1$ for $g$ given by
(\ref{eq:stabilizer-1}) and (\ref{eq:stabilizer-2}).
In the second case the differential is the identity map,
 in the first case
the differential is $dz\mapsto \alpha^*(dz)\alpha$.
We represent $\alpha$ as $p\Delta q$, where
$\Delta$ is diagonal with real eigenvalues and $p$, $q$ are unitary.
Now the statement becomes obvious.%
\hfill$\square$

\SS


{\bf\punct The degenerate principal series.
Proof of Proposition
\ref{pr:degenerate-principal}.
\label{ss:proof-degenerate}}
Thus, let $\Re(\sigma+\tau)=-n$, $\Im(\sigma)=\Im(\tau)=s$.
Then
$$
\det(a+uc)^{-n-\sigma|-n-\tau}=
|\det(a+uc)|^{-n-2is}
e^{i(\tau-\sigma)\mathrm{Arg}\det(a+uc)}
,$$
where $\mathrm{Arg}(\cdot)$
is the argument of a complex number.
Therefore
\begin{multline*}
\la T_{\sigma|\tau}(g)f_1,T_{\sigma|\tau}(g)f_2\ra_{L^2(\U(n))}
=\\=
\int_{\U(n)}f_1(u^{[g]})  \ov{f_2(u^{[g]}) }\,
\Bigl|\det(a+uc)^{-n-\sigma|-n-\tau}\Bigr|^2
\,d\mu(u)
=\\=
\int_{\U(n)} f_1(u^{[g]})  \ov{f_2(u^{[g]}) }
|\det(a+uc)|^{-2n}\,d\mu(u)
\end{multline*}
and we change the variable $z=u^{[g]}$ keeping in  mind
Lemma \ref{l:Jacobian}.%
\hfill$\square$

\SS


{\bf\punct The invariance of the  kernel.
Proof of Proposition \ref{pr:invariance-kernel}.
\label{ss:proof-invariance}}

\begin{lemma}
\label{l:identity-kernel}
  The distribution $\ell_{\sigma|\tau}$
 satisfies the identity
\begin{equation}
\ell_{\sigma|\tau} (u^{[g]}(v^{[g]})^*)
=\ell_{\sigma|\tau}(uv^*)
\det(a+u c)^{\{-\tau|-\sigma\}}
\det(a+vc)^{\{-\sigma|-\tau\}}
.
\label{eq:transformation-kernels}
\end{equation}
\end{lemma}

{\sc Proof.}
This follows from the identity
$$
1- u^{[g]}(v^{[g]})^*
=(a+uc)^{-1} (1-uv^*) (a+vc)^{*-1},\qquad
\text{where $g\in\U(n,n)$}
,$$
which can be easily verified by a straightforward
calculation (see e.g. \cite{Ner-notes}).
\hfill $\square$

\SS


{\sc Proof of Proposition \ref{pr:invariance-kernel}.
\label{ss:invariavce-kernel}}
 First, let $\Re(\sigma+\tau)<1$.
Substitute $h_1=u_1^{[g]}$, $h_2=u_2^{[g]}$ to the integral
$$
\la f_1,f_2\ra_{\sigma|\tau}=
\iint_{\U(n)\times\U(n)}
\ell_{\sigma|\tau}(h_1 h_2^*)\,f_1(h_1)\ov{f_2(h_2)}
\,d\mu(h_1)\,d\mu(h_2)
.$$
By the lemma, we obtain
\begin{multline*}
\iint_{\U(n)\times\U(n)}
\ell_{\sigma|\tau}(u_1u_2^*)
\det(a+u_1c)^{\{-\tau|-\sigma\}}
|\det(a+u_2c)|^{-\sigma|-\tau\}}
\times\\ \times
\,f_1(u_1)\ov{f_2(u_2)}
                            |\det(a+u_1c)|^{-2n}
                            |\det(a+u_2c)|^{-2n}
\,d\mu(u_1)\,d\mu(u_2)
=\\=
\la \rho_{\sigma|\tau}(g)f_1,
 \rho_{\sigma|\tau}(g)f_2\ra_{\sigma|\tau}
.\end{multline*}

Thus, our operators preserve the form
$\langle\cdot,\cdot\rangle_{\sigma|\tau}$.

 For general $\sigma$, $\tau\in\C$,
we consider the analytic continuation.
\hfill $\square$


\SS

{\bf\punct Shift of parameters.
 Proof of Proposition \ref{pr:shift}.
\label{ss:rroof-shift}}
First, we recall {\it Cartan decomposition}.
For $t_1\ge\dots\ge t_n$ denote
$$
\mathrm{CH}(t):
=
\begin{pmatrix}
\cosh(t_1)&0&\dots\\
0&\cosh(t_2)&\dots\\
\vdots&\vdots&\ddots
\end{pmatrix}
,
\quad
\mathrm{SH}(t):
=
\begin{pmatrix}
\sinh(t_1)&0&\dots\\
0&\sinh(t_2)&\dots\\
\vdots&\vdots&\ddots
\end{pmatrix}
.$$

The following statement is well known

\begin{proposition}
\label{pr:cartan}
Each element  $g\in\U(n,n)$ can be represented in the form
\begin{equation}
g=\begin{pmatrix}u_1&0\\0&v_1\end{pmatrix}
\begin{pmatrix} \mathrm{CH}(t)& \mathrm{SH}(t)\\
                 \mathrm{SH}(t)&\mathrm{CH}(t)
\end{pmatrix}
\begin{pmatrix}u_2&0\\0&v_2\end{pmatrix}
\label{eq:cartan}
\end{equation}
for some (uniquely determined) $t$ and some
$u_1$, $u_2$, $v_1$, $v_2\in\U(n)$.
\end{proposition}

Now we must show that the operator
$f(z)\mapsto \det(z)f(z)$ intertwines
$\rho_{\sigma|\tau}$ and $\rho_{\sigma+1|\tau-1}$
A straightforward calculation
reduces  this  to the identity
$$
\frac{\det(a+zc)}{\det\ov{(a+zc)}}
=\frac{\det(z^{[g]})}
{\det(z)}
,$$
which becomes obvious after the substitution
(\ref{eq:cartan}).

Also,
$$
\ell_{\sigma+1|\tau-1}(z)=
- \ell_{\sigma|\tau}(z) \det z
$$
and this easily implies the second statement of Proposition
\ref{pr:shift}.

\smallskip

\begin{figure}

\epsfbox{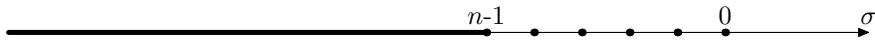}

\caption{Conditions of positivity of holomorphic
representations $\xi_\sigma$ (the '{\it Berezin--Wallach set}').}

\end{figure}

\section{Hilbert spaces of holomorphic functions}

\COUNTERS

 Theorem \ref{th:our-positivity}
exhaust the cases when the form
$\la\cdot,\cdot\ra_{\sigma|\tau}$
is positive definite on $C^\infty(\U(n))$.
However there are cases of positive semi-definiteness.
They are discussed in the next two sections.

Set $\tau=0$.
In this case,
 our construction produces holomorphic representations%
\footnote{or highest weight representations,}
of $\U(n,n)$. Holomorphic representations
were discovered by Harish-Chandra
(holomorphic discrete series, \cite{Harish15}) and Berezin
(analytic continuations of holomorphic discrete
series, \cite{Ber}).
They are discussed in numerous
 texts (for  partial expositions and further references,
see e.g. \cite{FK},
\cite{Ner-notes}), our aim is  to
show a link with our considerations.

\SS


{\bf\punct The case $\tau=0$.}
Substituting $\tau=0$, we get the action
$$
\rho_{\sigma|0}(g)\,f(z)=
f(z^{[g]})\det(a+zc)^{-n}\ov{\det(a+zc)}^{\,\,-n-\sigma}
.
$$
The Hermitian form is
$$
\la f_1,f_2\ra_{\sigma|0}=
\int_{\U(n)}\int_{\U(n)}
\det(1-z^* u)^\sigma f_1(z)\,\ov{f_2(u)}
\,d\mu(z)\,d\mu(u)
.
$$

\begin{theorem}
\label{th:ber-1}
The form $\la f_1,f_2\ra_{\sigma|0}$
 is positive semi-definite
iff $\sigma$ is contained in the set
$$
\text{$\sigma=0$, $-1$, \dots, $-(n-1)$,
or $\sigma<-(n-1)$}
.
$$
\end{theorem}

This means that all coefficients
$c_\bfm$ in the formula (\ref{eq:expansion-product})
are non-negative, but some coefficients vanish.
In fact the proof (see  below) is the examination
of these coefficients.

Under the conditions of the theorem we get a structure
of a pre-Hilbert space in $C^\infty(\U(n))$.
Denote by $\cH_\sigma$ the corresponding Hilbert
space.

Next, consider the action of the subgroup
$\U(n)\times\U(n)$ in $\cH_\sigma$.
We must get an orthogonal direct sum
$$
\bigoplus_{\bfm\in \Omega_\sigma} \pi_\bfm\oplus \pi_\bfm^*
$$
Some of summands of (\ref{eq:frobenius}) disappear,
when we pass to the quotient space; actually the summation
is taken over a proper subset $\Omega_\sigma$
of the set of all representations.
The next theorem
 is the description of
of the set $\Omega_\sigma$.

\begin{theorem}
a) If $\sigma<-(n-1)$, then
$$
\Omega_\sigma:=\Bigl\{\bfm:\, m_n\ge 0\Bigr\}
.
$$

b) If $\sigma=-n+\alpha$, where $\alpha=1$, $2$, \dots, $n-1$,
$n$, then
$$
\Omega_\sigma=\Bigl\{\bfm:m_n=0,\,m_{n-1}=1,\,\dots,\,
m_{n-\alpha+1}=\alpha-1\Bigr\}
.
$$
\end{theorem}

{\sc Proofs.}%
\footnote{This is the original Berezin's proof,
he started
from explicit expansions of reproducing kernels (\ref{eq:repro}).}
Substitute $\tau=0$ to (\ref{eq:expansion-2}),
\begin{multline}
c_{\bfm}=
(-1)^{n(n-1)/2}
2^{-\sigma n}
\prod_{j=1}^n\Gamma(\sigma+j)
\times
\\ \quad\times
\sum\limits_\bm\Biggl\{
\frac{(-1)^{\sum m_j}
\prod\limits_{1\le\alpha<\beta\le n} (m_\alpha-m_\beta)}
{\prod\limits_{j=1}^n\Gamma(\sigma-m_j+n)
   \boxed{ \Gamma(m_j+1)} }\chi_\bm(g) \Biggr\}
\label{eq:expansion-3}
=
\end{multline}
\begin{align}
=
\frac{(-1)^{n(n-1)/2}\sin^n(\pi\sigma)
2^{-(\sigma)n}}{\pi^n}
\prod_{j=1}^n\Gamma(\sigma+j)
\times
\qquad\qquad
\label{eq:expansion-31}
 \\
\qquad
\qquad\times
\sum\limits_\bm\Biggl\{
\prod\limits_{1\le\alpha<\beta\le n} (m_\alpha-m_\beta)
\prod\limits_{j=1}^n\frac{\Gamma(-\sigma+m_j-n+1)}
    {\boxed{\Gamma(m_j+1)}}\chi_\bm(g)
\Biggr\}
\label{eq:expansion-32}
\end{align}

We have
$\Gamma(m_j+1)=\infty$ for  $m_j<0$.
Therefore the corresponding fractions
in (\ref{eq:expansion-32}) are zero, and
the expansion of $\ell_{\sigma|0}$
has the form
\begin{equation}
\ell_{\sigma|0}=\sum_{\bfm:\, m_n\ge 0} c_{\bfm} \chi_{\bfm}
.
\label{eq:mn>0}
\end{equation}

Let us list possible cases.

\SS

{\sc Case 1.} If $\sigma<-n-1$, then all  coefficients
$c_{\bfm}$  are positive, see (\ref{eq:expansion-32});
in the line (\ref{eq:expansion-31}) poles of the Gamma-functions
cancel with zeros of sines.

\SS

{\sc Case 2.} If $\sigma\ge -n-1$ is non-integer,
 then all the coefficients
$c_\bfm$ are non-zero, but they have different signs.

\SS

\begin{figure}
a)

 \epsfbox{sobolev.15}

b)

 \epsfbox{sobolev.16}
\caption{ 'Maya diagrams' for signatures of harmonics
in holomorphic
representations.
\newline
a) A general case, $\sigma<n-1$.
\newline
b) Degenerate case. Here $\sigma=-(n-1)+5$.}
\end{figure}

{\sc Case 3.} Let  $\sigma$ be integer, $\sigma\ge -n+1$.
Consider a small perturbation of $\sigma$,
$$
\sigma=-n+\alpha+\epsilon
.
$$
In this case we get an uncertainty
in the expression (\ref{eq:expansion-3}):
$$
\frac{\prod_{j=1}^n\Gamma(-n+\alpha+\epsilon+j)}
{\prod_{j=1}^n\Gamma(\alpha-m_j+\epsilon)},
\qquad \epsilon\to 0
.$$
The order of the pole of the numerator is $n-\alpha$.
However order of a pole in the denominator
ranges between $n-\alpha$ and $n$ according
to $\bfm$. If the last
order $>n-\alpha$, then the ratio is zero.
The only possibility to get order of a pole $=n-\alpha$
is to set
\begin{equation}
m_n=0, \quad m_{n-1}=1,\quad\dots, \quad m_{n-\alpha+1}=0
\label{eq:signature-ber}
.\end{equation}
Thus the coefficients $c_{\bfm}$ are nonzero
only for signatures satisfying (\ref{eq:signature-ber});
they are positive.

We omit a discussion of positive integer $\sigma$
(the invariant inner product is not positive).
\hfill $\square$


\SS

{\bf\punct  Intertwining operators.}
Denote by $\B_n$ the space of complex $n\times n$-matrices
with norm $<1$.

Consider the integral operator
$$
I_\sigma f(z)=\int_{\U(n)} \det(1-zh^*)^\sigma f(h)\,
d\mu(h), \qquad z\in\B_n
.$$
It intertwines $\rho_{\sigma|0}$ with the representation
$\rho_{-n|-n-\sigma}$, denote the last representation
by $\xi_\sigma$
$$
\xi_\sigma(g) f(z)=f(z^{[g]}) \,\det (a+zc)^{\sigma}
.$$
The $I_\sigma$-image $\cH_\sigma^\circ$
 of the space $\cH_\sigma$ consists
of functions holomorphic in $\B_n$.
The structure of a Hilbert space in the space of holomorphic
functions is determined by the reproducing kernel
\begin{equation}
K_\alpha(z,\ov u)=\det(1-zu^*)^\sigma.
\label{eq:repro}
\end{equation}


{\bf\punct Concluding remarks (without proofs).}

\SS

a) For $\sigma<-(2n-1)$, the
 inner product in $\cH^\circ_\sigma$
can be written as an integral
$$
\la f_1, f_2\ra=\const \int_{\B_n}
 f_1(z)\,\ov{f_2(z)}\,\det(1-zz^*)^{-\sigma-2n}\,dz\,\ov{dz}
.$$

b) For $\sigma<n-1$ the space $\cH^\circ_\sigma$
contains all  polynomials.

\SS

c) Let $\sigma=0$, $-1$, \dots,$-(n-1)$.
Consider the matrix
$$
\Delta=
\begin{pmatrix}
\frac{\partial}{\partial z_{11}} & \dots &
                       \frac{\partial}{\partial z_{1n}}\\
\vdots&\ddots&\vdots \\
\frac{\partial}{\partial z_{n1}} & \dots &
                       \frac{\partial}{\partial z_{nn}}
\end{pmatrix}
.$$
The space $\cH^\circ_{-(n-1)}$ consists of functions $f$
satisfying the partial differential equation
$$
(\det \Delta)f(z)=0
.$$
The space $\cH^\circ_\sigma$, where $\sigma=0$, $-1$,
\dots, $-(n-1)$, consists of functions
that are annihilated by all
$(-\sigma+1)\times(-\sigma+1)$-minors of the matrix
$\Delta$. Also, $\cH_\sigma^\circ$ contains all  polynomial
satisfying this system of equations.

In particular, the space $\cH^\circ_0$ is one-dimensional.

\section{Unipotent representations.}

\COUNTERS

Here we propose models for 'unipotent representations'
of Sahi \cite{Sah4} and Dvorsky--Sahi,
\cite{DS1}--\cite{DS2}.


\smallskip

\begin{figure}

\epsfbox{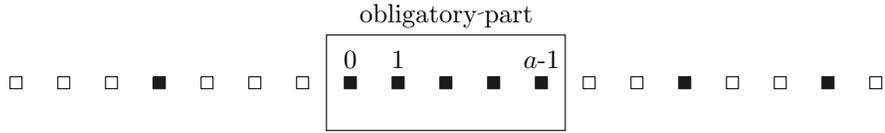}

\caption{Maya diagram for signatures $\in Z_j$;
here $j$ is the number of black boxes at the left
of the 'obligatory part'.}
\end{figure}

\begin{figure}

$$
a) \epsfbox{sobolev.9}
\qquad
b) \epsfbox{sobolev.10}
$$
$$
c) \epsfbox{sobolev.11}
\qquad
d) \epsfbox{sobolev.14}
$$

\caption{The case $n=2$.
\newline
a) $\alpha=0$. The decomposition of $L^2(\U(2)$ into
a direct sum.
\newline
b) $\alpha=1$. White circles correspond to the big
subrepresentation $W_{tail}$. The quotient is a direct sum of two
subrepresentations.
\newline
c) $\alpha=2$. The quotient is one-dimensional.
\newline
d) $0<\tau<1$, $\sigma=-n$. The invariant filtration.
The subquotients are unitary.}
\end{figure}


\SS

{\bf \punct Quotients of $\rho_{\sigma|\tau}$
at integer points.
\label{ss:unipotent}}
Set
\begin{equation}
\tau=0, \qquad \sigma=-n+\alpha,\qquad
\text{where $\alpha=0$, $1$, \dots, $n-1$.}
\end{equation}
For
$j=0$, $1$, \dots, $n-\alpha$ denote by
 $Z_j$ the set of all signatures $\bfm$
of the form
$$
\bfm=
 (m_1,\dots, m_{n-\alpha-j}, \alpha-1,\alpha-2,\dots,0,
m_{n-j+1},\dots m_n)
$$

Denote by $V_\bm$  the $\U(n)\times\U(n)$-subrepresentation
in $C^\infty(\U(n))$ corresponding a signature $\bm$,
see Subsection \ref{ss:compact}.

\smallskip

\begin{theorem}
\label{th:unipotent}
 The subspace
$$
W_{tail}:=\bigoplus_{\bm\notin \cup Z_j}
 V_\bm\subset C^\infty(\U(n)),
$$
is  $\U(n,n)$-invariant.

\smallskip

b)  The quotient $C^\infty(\U(n)) /W_{tail}$
is a sum $n-\alpha+1$ subrepresentations
$$
W_j=\oplus_{\bfm\in Z_j} V_ \bfm
.$$
The representation of $\U(n,n)$ in each
$W_j$ is unitary.
\end{theorem}

We formulate the result for $\alpha=0$
separately. In this case $W_{tail}=0$.

\begin{theorem}
\label{th:unipotent-2}
The representation $\rho_{-n|0}$ is a direct sum
of $n+1$ unitary representations $W_j$,
where $0\le j\le n$.
We have $V_\bfm\subset W_j$ if the number of
negative labels $m_k$ is $j$.
\end{theorem}

In particular, {\it we get a canonical decomposition of
$L^2(\U(n))$ into a direct sum of $(n+1)$ subspaces}.

\smallskip

Proof is given in the next subsection.

\smallskip

{\bf \punct The blow-up construction.
\label{ss:blow-2}}%
\footnote{The case $\U(1,1)$ was considered above
in Subsection \ref{ss:blow-1}.}
The distribution $\ell_{\sigma|\tau}$
depends meromorphically in
 two complex variables $\sigma$, $\tau$.
Its poles and zeros are located at $\sigma\in\Z$
and in $\tau\in\Z$. For this reason,
 values of $\ell_{\sigma|\tau}$
at points $(\sigma,\tau)\in\Z^2$ generally
 are not uniquely defined.
 Passing to such points from different directions,
 we get different limits%
\footnote{A remark for an expert  in algebraic
geometry: we consider blow up of the plane $\C^2$
at the point $(-n+\alpha,0)$.}.

Thus, set
\begin{equation}
\sigma=-n+\alpha+s\epsilon,\quad \tau=t\epsilon
\qquad \text{where $(s,t)\ne (0,0)$}.
\end{equation}
Substituting this to (\ref{eq:expansion-2}),
we get
\begin{multline}
\ell_{-n+\alpha+\epsilon s|\epsilon t}=
(-1)^{n(n-1)/2}
2^{-(\sigma+\tau)n}
\prod_{k=1}^n\Gamma\bigl(-n+\alpha \epsilon(s+t)+k\bigr)
\times
\\ \quad\times
\sum\limits_\bm\Biggl\{
\frac{(-1)^{\sum m_j}
\prod\limits_{1\le a<b\le n} (m_a-m_b)}
{\prod\limits_{k=1}^n\Gamma(\alpha+\epsilon s-m_k)
    \Gamma(\epsilon t+m_k+1)}\chi_\bm(g) \Biggr\}
\label{eq:expansion-4}
.\end{multline}

\begin{theorem}
\label{th:ell-s:t}
a) Let $s\ne-t$.
Then there exist a limit in the sense of distributions:
\begin{equation}
\ell^{s:t}(z):=
\lim_{\epsilon\to 0}
\ell_{-n+\alpha+\epsilon s|\epsilon t}(z)
,
\label{eq:ellst}
\end{equation}
In other words, the function
$(\sigma|\tau)\mapsto \ell_{\sigma|\tau}$
 has a removable singularity at
$\epsilon=0$ on the line
$$
\sigma=-n+\alpha+\epsilon s,
\qquad
\tau= \epsilon t,
\qquad \text{where $\epsilon \in\C$}
.
$$

b) Denote by $c_\bfm(s:t)$ the Fourier coefficients
of $\ell^{s:t}$.
 If $\bfm$ is in the 'tail', i.e., $\bfm\notin \cup Z_j$,
then $c_\bfm(s:t)=0$.

\SS

c) Moreover,
$\ell^{s:t}$ admits a decomposition
\begin{equation}
\ell^{s:t}=\sum_{j=0}^{n-\alpha} \frac{t^j s^{n-\alpha-j}}
{(s+t)^{n-\alpha}} \frL_j
\label{eq:s:t}
,
\end{equation}
where $\frL_j$ is of the form
\begin{equation}
\frL_j=\sum_{\bfm\in Z_j} a_\bfm \chi_\bfm
\label{eq:frL}
,\end{equation}
where $a_\bfm$ do not depend on $s$, $t$.

\SS

d) For each $j$ all  coefficients $a^j_\bfm$
in (\ref{eq:frL}) are either positive or negative.
\end{theorem}

{\sc Proof.}
 For the numerator of (\ref{eq:expansion-4})
we have the asymptotic
$$
\prod_{k=1}^n\Gamma\bigl(-n+\alpha \epsilon(s+t)+k\bigr)=
C\,\epsilon^{-n+\alpha}
(s+t)^{-n+\alpha}+O( \epsilon^{-n+\alpha+1}),
\qquad \epsilon\to 0
.
$$
Next, examine factors of the denominator,
$$
\Gamma(\alpha+\epsilon s-m_k)
    \Gamma(\epsilon t+m_k+1)\sim
\begin{cases}
A_1(m_k) (\epsilon t)^{-1} \quad&\text{if $m_k<0$}
\\
A_2(m_k)\quad &\text{if $0\le m_k<\alpha$}
\\
A_3(m_k) (\epsilon s)^{-1} &\text{if $m_k\ge \alpha$,}
\end{cases},\qquad\epsilon\to 0
$$
where $A_1$, $A_2$, $A_3$ do not depend on $s$, $t$.
Therefore, the order of the pole of denominator
$\prod_k$ of (\ref{eq:expansion-4})
 is
$$
\text{\it number of $m_j$ outside the segment $[0,\alpha-1]$}
$$
The minimal possible order of a pole of
the denominator is $n-\alpha$. In this case,
$c_\bfm$ has a finite nonzero limit,
of the form
$$
c_\bfm(s:t)= A(\bfm)\cdot\frac
{s^{\text{number of $m_k\ge\alpha$}}
\,\cdot\,
t^{\text{number of $m_k<0$.}} }
{(s+t)^{n-\alpha}}
$$

If an order of pole in the denominator is $> n-\alpha$, then
$c_\bfm(s:t)=0$. This corresponds to the tail.

We omit to watch the positivity of $c_\bfm(s:t)$.

Formally,  it is necessary to watch the growth
 of $c_\bfm(s:t)$ as $\bfm\to\infty$ and the growth of
$$
\frac{\partial}{\partial \epsilon}
c_\bfm(-n+\alpha+\epsilon s,\epsilon t)
$$
to be sure that (\ref{eq:ellst}) is a limit in the sense
of distributions. This is a more-or-less trivial exercise
on the Gamma-function.
\hfill $\square$

\smallskip

There are many ways to express $\frL^j$ in the terms
of $\ell^{s:t}$.  One of variants is given
in the following obvious proposition.

\begin{proposition}
The distribution $\frL_j$ is given by the formula
\begin{equation}
\frL_j(z)=\frac1 {j!}
\frac{\partial^j}{\partial t^j}
(1+t)^{n-\alpha} \ell^{1:t}(z)
\Bigr|_{t=0}
\end{equation}
\end{proposition}


{\bf \punct The family of invariant Hermitian forms.}
Thus, for $(\sigma, \tau)=(-n+\alpha,0)$
we obtained the following families of
$\rho_{-n+\alpha|0}$-invariant Hermitian
forms
\begin{equation}
R^{s:t}(f_1,f_2):=
\iint_{\U(n)\times\U(n)}
\ell^{s:t}(zu^*)f_1(z)\ov {f_2(u)}\,d\mu(z)\,d\mu(u)
\end{equation}
and
\begin{equation}
Q_j(f_1,f_2):=
\iint_{\U(n)\times\U(n)}
\frL_j(zu^*)f_1(z)\ov {f_2(u)}\,d\mu(z)\,d\mu(u)
\end{equation}

They are related as
$$
R^{s:t}(f_1,f_2)=\sum_{j=0}^{n-\alpha}\frac{t^j s^{n-\alpha-j}}
{(s+t)^{n-\alpha}} \frL_j(f_1,f_2)
$$

A form $\frL_j$ is zero on
$$
Y_j:=W_{tail}\oplus (\oplus_{i\ne j} W_i)
$$
and determines an inner product
on $W_j\simeq C^\infty(\U(n))/Y_j$.


\section{Some problems of harmonic analysis}

\COUNTERS

{\bf\punct Tensor products
 $\rho_{\sigma|\tau}\otimes\rho_{\sigma'|\tau'}$ .}
Nowadays the problem of decomposition of a tensor product
of two arbitrary unitary representations does not
seem interesting.
We propose several informal arguments
 for reasonableness
of the problem in our case.

\SS

a) For $n=1$ it is precisely the well-known
 problem of decomposition
of tensor products of unitary representations
of $\SL(2,\R)^\sim\simeq\SU(1,1)^\sim$,
 see \cite{Puk-tensor},
\cite{Mol}, \cite{Ros}, \cite{GKR}, \cite{Ner-jacobi}.

\SS

b) Decomposition of tensor
 products $\rho_{\sigma,0}\otimes\rho_{\sigma',0}$
of holomorphic representations
is a well-known combinatorial problem, see \cite{JV}.

\SS

c) Tensor products $\rho_{\sigma,0}\otimes\rho_{0|\tau'}$
are Berezin representations, see
\cite{Ber2}, \cite{UU}, \cite{Ner-berezin}.

\SS

d) All  problems a)--c) have interesting links
with theory of special functions.

\SS

e) There is a  canonical isomorphism
\footnote{Indeed, $\U(n)\simeq \U(n,n)/P$,
where $P$ is a maximal parabolic subgroup
in $\U(n,n)$.
The group $\U(n,n)$ has an open orbit on
$\U(n,n)/P\times\U(n,n)/P$, the  stabilizer
of a point is $\simeq\GL(n,\C)$.})
\begin{equation}
\rho_{-n/2|-n/2}\otimes \rho_{-n/2|-n/2}\simeq
L^2\bigl(\U(n,n)/\GL(n,\C)\bigr)
.
\end{equation}
Thus, we again come to a classical problem,
i.e., the problem of decomposition
of $L^2$ on a pseudo-Riemannian symmetric space $G/H$,
see \cite{FJ}, \cite{Osh}.%
\footnote{In a certain sense, the Plancherel formula
for $L^2(G/H)$ was obtained in \cite{vBS}, \cite{Del}.
However no Plancherel measure, nor spectra
are known. The corresponding problems remain open.}
General tensor products
$\rho_{\sigma|\tau}\otimes\rho_{\sigma'|\tau'}$
can be regarded as  deformations of the space
$L^2\bigl(\U(n,n)/\GL(n,\C)\bigr)$.


\SS

{\bf \punct Restriction problems.
\label{ss:restriction}}

\SS

1.Consider the group $G^*:=\U(n,n)$
and its subgroup $G:=\OO(n,n)$.
The group $G$ has an open dense orbit on
the space $\U(n)$, namely
$$
G/H:=\OO(n,n)/\OO(n,\C).
$$
The restriction of the representation
$\rho_{-n/2|-n/2}$ to $G$ is equivalent
to the representation
of $G$ in $L^2(G/H)$. Restrictions
of other $\rho_{\sigma|\tau}$
can be regarded as deformations of $L^2(G/H)$.

\SS

The same argument produces deformations
of $L^2$ on some other pseudo-Riemannian symmetric spaces.
Precisely, we have the following variants.

\SS

2. $G^*=\U(2n,2n)$, $G/H=\Sp(n,n)/\Sp(2n,\C)$.

\SS

3. $G^*=\U(n,n)$,  $G/H=\SOS(2n)/\OO(n,\C)$.

\SS

4. $G^*=\U(2n,2n)$, $G/H=\Sp(4n,\R)/\Sp(2n,\C)$.

\SS

5. $G^*=\U(p+q,p+q)$, $G/H=\U(p,q)\times\U(p,q)/\U(p,q)$.
In this case,  $G/H\simeq\U(p,q)$.

\SS

6. $G^*=\U(n,n)$, $G=\GL(n,\C)$.
In this case, we have $(n+1)$ open orbits
 $G/H_p=\GL(n,\C)/\U(p,n-p)$.

\SS

7. $G^*=\U(n,n)\times \U(n,n)$, $G=\U(n,n)$.
This is the problem about tensor products discussed
above.

\SS


{\bf\punct The Gelfand--Gindikin programm.}
Recall the statement of the problem,
see \cite{GG}, \cite{Ols}.
Let $G/H$ be a pseudo-Riemannian symmetric space.
The natural representation of $G$ in $L^2(G/H)$
has several pieces of spectrum.
Therefore, $L^2(G/H)$ admits a natural orthogonal decomposition
into direct summands having uniform spectra.
The problem is: {\it to
describe  explicitly the corresponding
subspaces or corresponding projectors.}

In Subsection \ref{ss:unipotent}
 we have obtained  a natural decomposition
of $L^2(\U(n))$ into $(n+1)$ direct summands.
Therefore {\it in the cases
listed in Subsection \ref{ss:restriction}
we have a natural orthogonal decompositions of $L^2(G/H)$}.

In any case, for the one-sheet hyperboloid
$\U(1,1)/\C^*$ we get the desired construction
(see Molchanov \cite{Mol-plane}, \cite{Ner-curve}).

\SS


{\bf\punct Matrix Sobolev spaces?}
Our inner product $\la\cdot,\cdot\ra_{\sigma|\tau}$
seems to be similar to Sobolev-type inner products
discussed in Subsection\ref{ss:sobolev-1}.
However it is not a Sobolev inner product, because
the kernel $\det(1-zu^*)^{\{\sigma|\tau\}}$
has a non-diagonal singularity.

 Denote
$$s=-\sigma-\tau+n.$$

 Let $F$ be a distribution on $\U(n)$,
let $F=\sum F_\bm$ be its expansion
 in a series of elementary harmonics.
We have
\begin{multline}
F\in \cH_{\sigma|\tau}
\qquad \Longleftrightarrow\qquad
\sum_\bm \frac{c_\bm}{\dim \pi_\bm} \|F_\bm\|^2_{L^2}<\infty
\qquad \Longleftrightarrow
\\  \Longleftrightarrow
\qquad
\sum_\bm\Bigl\{  \|F_\bm\|^2_{L^2}
\prod_{j=1}^n (1+|m_j|)^s\Bigr\}<\infty
\label{eq:matrix-sobolev}
,\end{multline}
where  $\|F_\bm\|_{L^2}$
 denotes
$$
\|F_\bm\|_{L^2}:=
\Bigl(\int_{\U(n)} |F_\bm(h)|^2\,d\mu(h)\Bigr)^{1/2}
$$
Our Hermitian form   defines a norm only in the case
$|s|<1$, but (\ref{eq:matrix-sobolev})
makes sense for  arbitrary real $s$.
Thus {\it we can  define
a Sobolev space ${\sf H}_s$ on $\U(n)$ of  arbitrary order}.

Author does not know applications
of this remark, but it seems that
 it can be useful in the  following
situation.

First, a reasonable harmonic
 analysis related to semisimple
Lie groups is the analysis
of unitary representations.
 But near 1980 Molchanov observed that
many identities with special function  admit
interpretations  on "physical level of rigor" as
formulas of non-unitary  harmonic analysis.
Up to now, there are no reasonable interpretations
of this phenomenon (but formulas exist,
see, e.g.
\cite{vDM}, see also
\cite{Ner-berezin}, Section 1.32 and formula (2.6)--(2.15) ).
In particular, we do not know
reasonable functional spaces that can be
scene  of action of such analysis.
 It seems that  our spaces
$\sf H_s$ can be possible candidates.

{\sf
University of Vienna, Math. Dept.,

Nordbergstarsse, 15,
Vienna,  Austria

\&
Math. Phys.Group, Institute for the Theoretical and Experimental Physics,

Bolshaya Cheremushkinskaya, 25,

Moscow 117 259, Russia

\& TFFA,  MechMath Dept., Moscow State University,

Vorob'evy Gory, Moscow, Russia

\tt neretin\@mccme.ru

URL:www.mat.univie.ac.at/$\sim$neretin,

    wwwth.itep.ru/$\sim$neretin

}


\begin{thebibliography}{cc}

\bibitem{vBS}
van den Ban, E.P., Schlichtkrull, H.,
{\it The most continuous part of the Plancherel decomposition
for a reductive symmetric space.}
Ann. Math., 145 (1997), 267--364

\bibitem{Bar}
Bargmann, V.
{\it Irreducible unitary representations of
the Lorentz group.} Ann. Math, 48 (1947), 568--640

\bibitem{Ber}
Berezin, F.A., {\it Quantization
in complex symmetric spaces.}
Izv. Akad. Nauk SSSR, Ser. Math., 39, 2, 1362--1402 (1975);
English translation: Math USSR Izv. 9 (1976),
 No 2, 341--379(1976)

\bibitem{Ber2}
Berezin, F. A. {\it The connection between covariant
 and contravariant symbols of operators
on classical complex symmetric spaces.}
 Sov. Math. Dokl. 19 (1978), 786--789


\bibitem{BOO}
Branson, Th., Olafsson, G., Orsted, B.
{\it Spectrum generating operators and intertwining operators
for representations induced from a maximal parabolic subgroup,}
J. Funct. Anal., 135, 163--205.

\bibitem{Del}
Delorme, P.
{\it Formule de Plancherel pour les espaces symm\'etrique
reductifs.} Ann. Math., 147 (1998), 417--452



\bibitem{vDM}
van Dijk, G., Molchanov, V.F.
{\it The Berezin form for rank
one para-Hermitian symmetric spaces.}
J. Math. Pure. Appl., 78 (1999),  99--119.

\bibitem{DS1}
  Dvorsky, A., Sahi, S.
 {\it  Explicit Hilbert spaces for certain
unipotent representations. II.}
Invent. Math. 138 (1999), no. 1, 203--224.




\bibitem{DS2}
 Dvorsky, A., Sahi, S.
 {\it Explicit Hilbert spaces
  for certain unipotent representations. III.}
 J. Funct. Anal.
   201(2003),  no. 2, 430--456.

\bibitem{FK}
Faraut, J., Koranyi, A.,
 {\it Analysis in symmetric cones.}
Oxford Univ.Press, (1994)

\bibitem{FJ}
 Flensted-Jensen, M.
 {\it Discrete series for semisimple symmetric spaces.} Ann.
of Math. (2) 111 (1980), no. 2, 253--311.

\bibitem{Fri}
Friedrichs, K. O. {\it  Mathematical aspects
 of the quantum theory of fields.}
Interscience Publishers, Inc., New York,
 1953.

\bibitem{GG}
Gelfand, I. M., Gindikin, S. G.
{\it Complex manifolds whose skeletons
are semisimple Lie groups and analytic discrete series
of representations.}
Funct. Anal. Appl., 11 (1978), 258--265


\bibitem{GN2}
 Gelfand, I.M., Naimark, M.I.,
{\it Unitary representations of
classical groups.}
 {\it Unitary representations of classical groups.}
 Trudy MIAN., t.36 (1950);
 German translation: Gelfand I.N., Neumark M.A.,
 {\it Unitare Darstellungen der klassischen gruppen.},
 Akademie-Verlag, Berlin, 1957.




\bibitem{GR}
Gradshtein, I.S., Ryzhik, I.M.
{\it Tables of integrals, sums and products.}
Fizmatgiz, 1963; English translation: Acad. Press, NY, 1965

\bibitem{GKR}
Groenevelt, W., Koelink, E., Rosengren, H.
 {\it Continuous Hahn polynomials
and Clebsch--Gordan coefficients.}
Preprint {\tt http://arxiv.org/math/abs/0302251}

\bibitem{Harish15}
 Harish-Chandra, {\it Representations of semisimple Lie
groups IV,} Amer. J. Math.,  743--777 (1955).
Reprinted in Harish-Chandra {\it Collected papers,}
v.2.

\bibitem{HTF}
{\it Higher transcendental functions}, v.1.,
McGraw-Hill book company, 1953




\bibitem{JV}
Jakobsen, H.P., Vergne, M.,
{\it Restrictions and expansions of holomorphic
representations.}
J. Funct. Anal., 34 (1979), 29--53.

\bibitem{Kad}
Kadell, K. {\it The Selberg--Jack symmetric functions.}
Adv. Math., 130 (1997), 33-102



\bibitem{Kir}
Kirillov, A.A.
{\it Elements of representation theory,}
Moscow, Nauka, 1972; English transl.:
Springer, 1976.


\bibitem{Kra}
Krattenthaler, C.
{\it Advanced determinant calculus.}
The Andrews Festschrift (Maratea, 1998).
Sem. Lothar. Combin. 42 (1999), Art. B42q, 67 pp.
 (electronic).


\bibitem{Mol}
 Molchanov, V. F.{\it Tensor products
 of unitary representations of the
three-dimensional Lorentz group.}
 Izv. Akad. Nauk SSSR Ser. Mat. 43 (1979), no. 4,
860--891, 967. English transl. in Izvestia.

\bibitem{Mol-plane}
 Molchanov, V. F.
{\it Quantization on the imaginary Lobachevsky plane.}
Funct. Anal. Appl., 14 (1980), 162--144





\bibitem{Ner-curve}
Neretin, Yu. A. {\it The restriction of functions
holomorphic in a domain to a curve lying in the boundary,
and discrete $\SL_2(\R)$-spectra.}
Izvestia: Mathematics, 62:3(1998), 493--513

\bibitem{Ner-beta}      
Neretin, Yu.A.,
 {\it Matrix analogs of $\B$-function
 and Plancherel formula for Berezin kernel representations,}
Mat. Sbornik, 191, No.5 (2000), 67--100;

\bibitem{Ner-berezin}     
Neretin, Yu.A.,
{\it Plancherel formula for Berezin deformation
 of $L^2$ on Riemannian
symmetric space},
 J. Funct. Anal. (2002), 189(2002), 336--408.

\bibitem{Ner-determinant}
Neretin, Yu.A.
{\it Matrix balls, radial analysis
 of Berezin kernels,
and hypergeometric determinants,}
Moscow Math. J., v.1 (2001), 157--221.


\bibitem{Ner-stein}
Neretin, Yu.A.
{\it Notes Sahi--Stein representations
and some problems of non-$L^2$ harmonic analysis.},
J. Math. Sci., New York, 141 (2007), 1452--1478


\bibitem{Ner-preprint}
Neretin, Yu. A. {\it Notes on matrix analogs
 of Sobolev spaces
and Stein--Sahi representations.}
Preprint, {\tt http://arxiv.org/abs/math/0411419}

\bibitem{Ner-notes}
Neretin, Yu. A. {\it Lectures on Gaussian
integral operators and classical groups,} to appear.

\bibitem{Ner-jacobi}
Neretin, Yu. A.
{\it Some continuous analogs of expansion in Jacobi
polynomials and vector-valued orthogonal bases.}
Funct. Anal. Appl., 39 (2005), 31--46.

\bibitem{NO}
Neretin, Yu.A., Olshanskii, G.I.,
{\it  Boundary values of
 holomorphic functions ,
 singular unitary representations of groups
$O(p,q)$  and their limits as  $q\to\infty$.}
Zapiski nauchn. semin. POMI RAN 223, 9--91(1995);
 English translation:
J.Math.Sci., New York, 87, 6 (1997), 3983--4035.


\bibitem{Ols}
Olshanskij, G.I., {\it Complex Lie semigroups, Hardy
spaces, and Gelfand--Gindikin programm.}
Deff. Geom. Appl., 1 (1991), 235--246

\bibitem{Osh}
Oshima, T.
{\it A calculation of $c$-functions
 for semisimple symmetric spaces.
Lie groups and symmetric spaces,}
  307--330, Amer. Math. Soc. Transl. Ser. 2, 210,
 Amer. Math. Soc.,
   Providence, RI, 2003.

\bibitem{Puk-tensor}
Pukanszky, L., {\it On the Kronecker products of irreducible
 unitary representations  of the
 $2\times2$ real unimodular  group.}
  Trans. Amer. Math. Soc., 100
(1961), 116--152

\bibitem{Puk-classic}
Pukanzsky, L. {\it  Plancherel formula for
universal covering group of $\SL(2,\R)$.}
Math. Ann., 156 (1964), 96-?

 \bibitem{RS}
  Ricci, F., Stein, E. M.
  {\it Homogeneous distributions
   on spaces of Hermitean matrices.}
    J. Reine Angew. Math.  368  (1986), 142--164.

\bibitem{Ros} Rosengren, H.
{\it Multilinear Hankel forms of higher order
and orthogonal polynomials.} Math. Scand., 82
(1998), 53-88

\bibitem{Sah1}
   Sahi, S.
  {\it  A simple construction of Stein's complementary series
representations.} Proc. Amer. Math. Soc. 108 (1990),
 no. 1, 257--266.


\bibitem{Sah2}
 Sahi, S.,
 {\it Unitary representations on the Shilov
 boundary of a symmetric tube domain,}
  Contemp. Math. 145 (1993) 275--286.

\bibitem{Sah3}
 Sahi, S. {\it Jordan algebras
 and degenerate principal series}, J. Reine
 Angew.Math. 462 (1995) 1--18.

 \bibitem{Sah4}
 Sahi, S. {\it Explicit Hilbert spaces
 for certain unipotent representations.}
   Invent. Math.  110  (1992),  no. 2, 409--418.

\bibitem{SS}
 Sahi, S.,  Stein, E. M.
{\it Analysis in matrix space and Speh's
representations.} Invent. Math. 101 (1990), no. 2, 379--393.


\bibitem{Sal}
Sally, P. J., {\it Analytic continuations of irreducible
unitary representations of
the universal covering group of $\SL(2,\R)$.}
Amer. Math. Soc., Providence, 1967

\bibitem{Ste}
 Stein, E. M.
 {\it Analysis in matrix spaces and
  some new representations of ${\rm SL}(N,\,C)$.}
 Ann. of Math. (2) 86 1967 461--490.



\bibitem{UU}
Unterberger, A., Upmeier, H., {\it The Berezin transform and
 invariant differential operators}.
Comm.Math.Phys.,164, 563--597(1994)

\bibitem{Vog}
Vogan, D. A., {\it The unitary dual
 of ${\rm GL}(n)$ over an Archimedean field.}
   Invent. Math.  83  (1986),  no. 3, 449--505.


\bibitem{Wey}
 Weyl, H. {\it The Classical Groups.
 Their Invariants and Representations.} Princeton
University Press, Princeton, N.J., 1939.


\end{thebibliography}
\end{document}